\documentclass{gtart}

\def\ifplaintex{\expandafter\ifx\csname documentclass\endcsname\relax}

\def\gtp{{\mathsurround=0pt\it $\cal G\mskip-2mu$eometry \&\ 
$\cal T\!\!$opology $\cal P\!$ublications}}  

\def\recd{{\small Received:\qua\receiveddate\ifx\reviseddate\relax
\else\qquad Revised:\qua\reviseddate\fi\par}} 

\def\lognumber#1{\def\thelognumber{#1}}
\def\volumenumber#1{\def\thevolumenumber{#1}}
\def\volumeyear#1{\def\thevolumeyear{#1}}
\def\papernumber#1{\def\thepapernumber{#1}}
\def\pagenumbers#1#2{\def\startpage{#1}\def\finishpage{#2}}
\def\published#1{\def\publishdate{#1}}

\def\received#1{\def\receiveddate{#1}}

\def\accepted#1{\def\accepteddate{#1}}

\long\def\asciiabstract#1{\long\def\theasciiabstract{#1}}

\let\\\par\let\thelognumber\relax\let\thevolumenumber\relax
\let\thepapernumber\relax\let\thevolumeyear\relax\let\startpage\relax
\let\finishpage\relax\let\publishdate\relax\let\receiveddate\relax
\let\reviseddate\relax\let\accepteddate\relax\let\theasciititle\relax
\let\theasciiauthors\relax
\let\theasciiabstract\relax

\let\theasciiemail\relax

\ifplaintex
\font\logobig=cmssbx10 scaled 3836
\font\logomed=cmssbx10 scaled 2557
\else
\font\logobig=cmssbx10 scaled 4200
\font\logomed=cmssbx10 scaled 2800
\fi

\long\def\makeagttitle{   
\count0=\startpage
\agt\hfill      
\hbox to 45truept{\vbox to 0pt{\vglue -13truept{\logomed A\kern -.37em{\logobig 
T}\kern -.38em G}\vss}\hss}
\break
{\small Volume \thevolumenumber\ (\thevolumeyear)
\startpage--\finishpage\nl
Published: \publishdate}

\vglue .25truein

{\parskip=0pt\leftskip 0pt plus
1fil\def\\{\par\smallskip}{\Large\bf\thetitle}\par\medskip} \vglue
0.05truein

{\parskip=0pt\leftskip 0pt plus 1fil\def\\{\par}{\sc\theauthors}
\par\medskip}%
 
\vglue 0.03truein 

{\small\leftskip 25truept\rightskip 25truept{\bf Abstract}\stdspace\theabstract

{\bf AMS Classification}\stdspace\theprimaryclass
\ifx\thesecondaryclass\relax\else; \thesecondaryclass\fi\par
{\bf Keywords}\stdspace \thekeywords\par}\vglue 7truept

}   

\ifplaintex
\font\phead=cmsl9 scaled 950
\font\pnum=cmbx10 scaled 913
\font\pfoot=cmsl9 scaled 950
\headline{\vbox to 0pt{\vskip -4.5mm\line{\small\phead\ifnum
\count0=\startpage ISSN 1472-2739 (on-line) 1472-2747 (printed)
\hfill {\pnum\folio}\else\ifodd\count0\def\\{ }%
\ifx\theshorttitle\relax\thetitle\else\theshorttitle\fi\hfill{\pnum\folio}
\else\def\\{ and }{\pnum\folio}\hfill\ifx\theshortauthors\relax\theauthors
\else\theshortauthors\fi\fi\fi}\vss}}
\footline{\vbox to 0pt{\vglue 0mm\line{\small\pfoot\ifnum\count0=\startpage
\copyright\ \gtp\hfill\else
\agt, Volume \thevolumenumber\ (\thevolumeyear)\hfill\fi}\vss}}
\else
\font=cmsl9 scaled 1050
\font\lnum=cmbx10 
\font=cmsl9 scaled 1050
\makeatletter
\def\@oddhead{{\small\lhead\ifnum\count0=\startpage ISSN 1472-2739 
(on-line) 1472-2747 (printed)\hfill {\lnum\number\count0}\else\ifodd\count0
\def\\{ }\ifx\theshorttitle\relax \thetitle \else\theshorttitle\fi\hfill
{\lnum\number\count0}\else\def\\{ and }{\lnum\number\count0}
\hfill\ifx\theshortauthors\relax 
\theauthors\else\theshortauthors\fi\fi\fi}}\def\@evenhead{\@oddhead}
\def\@oddfoot{\small\lfoot\ifnum\count0=\startpage\copyright\ \gtp\hfill\else
\agt, Volume \thevolumenumber\ (\thevolumeyear)\hfill\fi}
\def\@evenfoot{\@oddfoot}
\makeatother
\fi
\let\maketitlepage\makeagttitle
\let\makeshorttitle\maketitlepage
\let\maketitle\maketitlepage

\newwrite\gtoutfile
\long\gdef\makeheadfile{  
{\def\\{, }\def\s{ }
\immediate\openout\gtoutfile head.xxx
\immediate\write\gtoutfile{To: math@arxiv.org}
\immediate\write\gtoutfile{Subject: put OR rep NNNNN:ppppp}
\immediate\write\gtoutfile{--text follows this line--}
\immediate\write\gtoutfile{Proxy-for: \ifx\theasciiauthors\relax
\theauthors\else\theasciiauthors\fi\s<\ifx\theasciiemail\relax\theemail\else\theasciiemail\fi>}
\immediate\write\gtoutfile{\noexpand\\}
\immediate\write\gtoutfile{Authors: \ifx\theasciiauthors\relax
\theauthors\else\theasciiauthors\fi}
{\def\\{ }\immediate\write\gtoutfile{Title: \ifx\theasciititle\relax
\thetitle\else\theasciititle\fi}}
\immediate\write\gtoutfile{Subj-class: GT or SG, GR etc}
\immediate\write\gtoutfile{MSC-class: \theprimaryclass\ifx\thesecondaryclass\relax\else, \thesecondaryclass\fi}
\immediate\write\gtoutfile{Journal-ref: Algebr. Geom. Topol. \thevolumenumber\s
(\thevolumeyear) \startpage-\finishpage}
\immediate\write\gtoutfile{Comments: Published by Algebraic and
Geometric Topology at}
\immediate\write\gtoutfile{\s\s\s  http://www.maths.warwick.ac.uk/agt/AGTVol\thevolumenumber/agt-\thevolumenumber-\thepapernumber.abs.html}
\immediate\write\gtoutfile{\noexpand\\}
\immediate\write\gtoutfile{}
\ifx\theasciiabstract\relax
\immediate\write\gtoutfile{\theabstract}\else
\immediate\write\gtoutfile{\theasciiabstract}\fi
\immediate\write\gtoutfile{}
\immediate\write\gtoutfile{\noexpand\\}
\immediate\write\gtoutfile{}
\immediate\closeout\gtoutfile}}  

\def\maketitlepage{\makeagttitle\makeheadfile}
\let\makeshorttitle\maketitlepage
\let\maketitle\maketitlepage

\def\ifplaintex{\expandafter\ifx\csname documentclass\endcsname\relax}

\def\gtp{{\mathsurround=0pt\it $\cal G\mskip-2mu$eometry \&\ 
$\cal T\!\!$opology $\cal P\!$ublications}}  

\def\recd{{\small Received:\qua\receiveddate\ifx\reviseddate\relax
\else\qquad Revised:\qua\reviseddate\fi\par}} 

\def\lognumber#1{\def\thelognumber{#1}}
\def\volumenumber#1{\def\thevolumenumber{#1}}
\def\volumeyear#1{\def\thevolumeyear{#1}}
\def\papernumber#1{\def\thepapernumber{#1}}
\def\pagenumbers#1#2{\def\startpage{#1}\def\finishpage{#2}}
\def\published#1{\def\publishdate{#1}}

\def\received#1{\def\receiveddate{#1}}

\def\accepted#1{\def\accepteddate{#1}}

\long\def\asciiabstract#1{\long\def\theasciiabstract{#1}}

\let\\\par\let\thelognumber\relax\let\thevolumenumber\relax
\let\thepapernumber\relax\let\thevolumeyear\relax\let\startpage\relax
\let\finishpage\relax\let\publishdate\relax\let\receiveddate\relax
\let\reviseddate\relax\let\accepteddate\relax\let\theasciititle\relax
\let\theasciiauthors\relax
\let\theasciiabstract\relax

\let\theasciiemail\relax

\ifplaintex
\font\logobig=cmssbx10 scaled 3836
\font\logomed=cmssbx10 scaled 2557
\else
\font\logobig=cmssbx10 scaled 4200
\font\logomed=cmssbx10 scaled 2800
\fi

\long\def\makeagttitle{   
\count0=\startpage
\agt\hfill      
\hbox to 45truept{\vbox to 0pt{\vglue -13truept{\logomed A\kern -.37em{\logobig 
T}\kern -.38em G}\vss}\hss}
\break
{\small Volume \thevolumenumber\ (\thevolumeyear)
\startpage--\finishpage\nl
Published: \publishdate}

\vglue .25truein

{\parskip=0pt\leftskip 0pt plus
1fil\def\\{\par\smallskip}{\Large\bf\thetitle}\par\medskip} \vglue
0.05truein

{\parskip=0pt\leftskip 0pt plus 1fil\def\\{\par}{\sc\theauthors}
\par\medskip}%
 
\vglue 0.03truein 

{\small\leftskip 25truept\rightskip 25truept{\bf Abstract}\stdspace\theabstract

{\bf AMS Classification}\stdspace\theprimaryclass
\ifx\thesecondaryclass\relax\else; \thesecondaryclass\fi\par
{\bf Keywords}\stdspace \thekeywords\par}\vglue 7truept

}   

\ifplaintex
\font\phead=cmsl9 scaled 950
\font\pnum=cmbx10 scaled 913
\font\pfoot=cmsl9 scaled 950
\headline{\vbox to 0pt{\vskip -4.5mm\line{\small\phead\ifnum
\count0=\startpage ISSN 1472-2739 (on-line) 1472-2747 (printed)
\hfill {\pnum\folio}\else\ifodd\count0\def\\{ }%
\ifx\theshorttitle\relax\thetitle\else\theshorttitle\fi\hfill{\pnum\folio}
\else\def\\{ and }{\pnum\folio}\hfill\ifx\theshortauthors\relax\theauthors
\else\theshortauthors\fi\fi\fi}\vss}}
\footline{\vbox to 0pt{\vglue 0mm\line{\small\pfoot\ifnum\count0=\startpage
\copyright\ \gtp\hfill\else
\agt, Volume \thevolumenumber\ (\thevolumeyear)\hfill\fi}\vss}}
\else
\font=cmsl9 scaled 1050
\font\lnum=cmbx10 
\font=cmsl9 scaled 1050
\makeatletter
\def\@oddhead{{\small\lhead\ifnum\count0=\startpage ISSN 1472-2739 
(on-line) 1472-2747 (printed)\hfill {\lnum\number\count0}\else\ifodd\count0
\def\\{ }\ifx\theshorttitle\relax \thetitle \else\theshorttitle\fi\hfill
{\lnum\number\count0}\else\def\\{ and }{\lnum\number\count0}
\hfill\ifx\theshortauthors\relax 
\theauthors\else\theshortauthors\fi\fi\fi}}\def\@evenhead{\@oddhead}
\def\@oddfoot{\small\lfoot\ifnum\count0=\startpage\copyright\ \gtp\hfill\else
\agt, Volume \thevolumenumber\ (\thevolumeyear)\hfill\fi}
\def\@evenfoot{\@oddfoot}
\makeatother
\fi
\let\maketitlepage\makeagttitle
\let\makeshorttitle\maketitlepage
\let\maketitle\maketitlepage

\newwrite\gtoutfile
\long\gdef\makeheadfile{  
{\def\\{, }\def\s{ }
\immediate\openout\gtoutfile head.xxx
\immediate\write\gtoutfile{To: math@arxiv.org}
\immediate\write\gtoutfile{Subject: put OR rep NNNNN:ppppp}
\immediate\write\gtoutfile{--text follows this line--}
\immediate\write\gtoutfile{Proxy-for: \ifx\theasciiauthors\relax
\theauthors\else\theasciiauthors\fi\s<\ifx\theasciiemail\relax\theemail\else\theasciiemail\fi>}
\immediate\write\gtoutfile{\noexpand\\}
\immediate\write\gtoutfile{Authors: \ifx\theasciiauthors\relax
\theauthors\else\theasciiauthors\fi}
{\def\\{ }\immediate\write\gtoutfile{Title: \ifx\theasciititle\relax
\thetitle\else\theasciititle\fi}}
\immediate\write\gtoutfile{Subj-class: GT or SG, GR etc}
\immediate\write\gtoutfile{MSC-class: \theprimaryclass\ifx\thesecondaryclass\relax\else, \thesecondaryclass\fi}
\immediate\write\gtoutfile{Journal-ref: Algebr. Geom. Topol. \thevolumenumber\s
(\thevolumeyear) \startpage-\finishpage}
\immediate\write\gtoutfile{Comments: Published by Algebraic and
Geometric Topology at}
\immediate\write\gtoutfile{\s\s\s  http://www.maths.warwick.ac.uk/agt/AGTVol\thevolumenumber/agt-\thevolumenumber-\thepapernumber.abs.html}
\immediate\write\gtoutfile{\noexpand\\}
\immediate\write\gtoutfile{}
\ifx\theasciiabstract\relax
\immediate\write\gtoutfile{\theabstract}\else
\immediate\write\gtoutfile{\theasciiabstract}\fi
\immediate\write\gtoutfile{}
\immediate\write\gtoutfile{\noexpand\\}
\immediate\write\gtoutfile{}
\immediate\closeout\gtoutfile}}  

\def\maketitlepage{\makeagttitle\makeheadfile}
\let\makeshorttitle\maketitlepage
\let\maketitle\maketitlepage

\def\ifplaintex{\expandafter\ifx\csname documentclass\endcsname\relax}

\def\gtp{{\mathsurround=0pt\it $\cal G\mskip-2mu$eometry \&\ 
$\cal T\!\!$opology $\cal P\!$ublications}}  

\def\recd{{\small Received:\qua\receiveddate\ifx\reviseddate\relax
\else\qquad Revised:\qua\reviseddate\fi\par}} 

\def\lognumber#1{\def\thelognumber{#1}}
\def\volumenumber#1{\def\thevolumenumber{#1}}
\def\volumeyear#1{\def\thevolumeyear{#1}}
\def\papernumber#1{\def\thepapernumber{#1}}
\def\pagenumbers#1#2{\def\startpage{#1}\def\finishpage{#2}}
\def\published#1{\def\publishdate{#1}}

\def\received#1{\def\receiveddate{#1}}

\def\accepted#1{\def\accepteddate{#1}}

\long\def\asciiabstract#1{\long\def\theasciiabstract{#1}}

\let\\\par\let\thelognumber\relax\let\thevolumenumber\relax
\let\thepapernumber\relax\let\thevolumeyear\relax\let\startpage\relax
\let\finishpage\relax\let\publishdate\relax\let\receiveddate\relax
\let\reviseddate\relax\let\accepteddate\relax\let\theasciititle\relax
\let\theasciiauthors\relax
\let\theasciiabstract\relax

\let\theasciiemail\relax

\ifplaintex
\font\logobig=cmssbx10 scaled 3836
\font\logomed=cmssbx10 scaled 2557
\else
\font\logobig=cmssbx10 scaled 4200
\font\logomed=cmssbx10 scaled 2800
\fi

\long\def\makeagttitle{   
\count0=\startpage
\agt\hfill      
\hbox to 45truept{\vbox to 0pt{\vglue -13truept{\logomed A\kern -.37em{\logobig 
T}\kern -.38em G}\vss}\hss}
\break
{\small Volume \thevolumenumber\ (\thevolumeyear)
\startpage--\finishpage\nl
Published: \publishdate}

\vglue .25truein

{\parskip=0pt\leftskip 0pt plus
1fil\def\\{\par\smallskip}{\Large\bf\thetitle}\par\medskip} \vglue
0.05truein

{\parskip=0pt\leftskip 0pt plus 1fil\def\\{\par}{\sc\theauthors}
\par\medskip}%
 
\vglue 0.03truein 

{\small\leftskip 25truept\rightskip 25truept{\bf Abstract}\stdspace\theabstract

{\bf AMS Classification}\stdspace\theprimaryclass
\ifx\thesecondaryclass\relax\else; \thesecondaryclass\fi\par
{\bf Keywords}\stdspace \thekeywords\par}\vglue 7truept

}   

\ifplaintex
\font\phead=cmsl9 scaled 950
\font\pnum=cmbx10 scaled 913
\font\pfoot=cmsl9 scaled 950
\headline{\vbox to 0pt{\vskip -4.5mm\line{\small\phead\ifnum
\count0=\startpage ISSN 1472-2739 (on-line) 1472-2747 (printed)
\hfill {\pnum\folio}\else\ifodd\count0\def\\{ }%
\ifx\theshorttitle\relax\thetitle\else\theshorttitle\fi\hfill{\pnum\folio}
\else\def\\{ and }{\pnum\folio}\hfill\ifx\theshortauthors\relax\theauthors
\else\theshortauthors\fi\fi\fi}\vss}}
\footline{\vbox to 0pt{\vglue 0mm\line{\small\pfoot\ifnum\count0=\startpage
\copyright\ \gtp\hfill\else
\agt, Volume \thevolumenumber\ (\thevolumeyear)\hfill\fi}\vss}}
\else
\font=cmsl9 scaled 1050
\font\lnum=cmbx10 
\font=cmsl9 scaled 1050
\makeatletter
\def\@oddhead{{\small\lhead\ifnum\count0=\startpage ISSN 1472-2739 
(on-line) 1472-2747 (printed)\hfill {\lnum\number\count0}\else\ifodd\count0
\def\\{ }\ifx\theshorttitle\relax \thetitle \else\theshorttitle\fi\hfill
{\lnum\number\count0}\else\def\\{ and }{\lnum\number\count0}
\hfill\ifx\theshortauthors\relax 
\theauthors\else\theshortauthors\fi\fi\fi}}\def\@evenhead{\@oddhead}
\def\@oddfoot{\small\lfoot\ifnum\count0=\startpage\copyright\ \gtp\hfill\else
\agt, Volume \thevolumenumber\ (\thevolumeyear)\hfill\fi}
\def\@evenfoot{\@oddfoot}
\makeatother
\fi
\let\maketitlepage\makeagttitle
\let\makeshorttitle\maketitlepage
\let\maketitle\maketitlepage

\newwrite\gtoutfile
\long\gdef\makeheadfile{  
{\def\\{, }\def\s{ }
\immediate\openout\gtoutfile head.xxx
\immediate\write\gtoutfile{To: math@arxiv.org}
\immediate\write\gtoutfile{Subject: put OR rep NNNNN:ppppp}
\immediate\write\gtoutfile{--text follows this line--}
\immediate\write\gtoutfile{Proxy-for: \ifx\theasciiauthors\relax
\theauthors\else\theasciiauthors\fi\s<\ifx\theasciiemail\relax\theemail\else\theasciiemail\fi>}
\immediate\write\gtoutfile{\noexpand\\}
\immediate\write\gtoutfile{Authors: \ifx\theasciiauthors\relax
\theauthors\else\theasciiauthors\fi}
{\def\\{ }\immediate\write\gtoutfile{Title: \ifx\theasciititle\relax
\thetitle\else\theasciititle\fi}}
\immediate\write\gtoutfile{Subj-class: GT or SG, GR etc}
\immediate\write\gtoutfile{MSC-class: \theprimaryclass\ifx\thesecondaryclass\relax\else, \thesecondaryclass\fi}
\immediate\write\gtoutfile{Journal-ref: Algebr. Geom. Topol. \thevolumenumber\s
(\thevolumeyear) \startpage-\finishpage}
\immediate\write\gtoutfile{Comments: Published by Algebraic and
Geometric Topology at}
\immediate\write\gtoutfile{\s\s\s  http://www.maths.warwick.ac.uk/agt/AGTVol\thevolumenumber/agt-\thevolumenumber-\thepapernumber.abs.html}
\immediate\write\gtoutfile{\noexpand\\}
\immediate\write\gtoutfile{}
\ifx\theasciiabstract\relax
\immediate\write\gtoutfile{\theabstract}\else
\immediate\write\gtoutfile{\theasciiabstract}\fi
\immediate\write\gtoutfile{}
\immediate\write\gtoutfile{\noexpand\\}
\immediate\write\gtoutfile{}
\immediate\closeout\gtoutfile}}  

\def\maketitlepage{\makeagttitle\makeheadfile}
\let\makeshorttitle\maketitlepage
\let\maketitle\maketitlepage

\def\ifplaintex{\expandafter\ifx\csname documentclass\endcsname\relax}

\def\gtp{{\mathsurround=0pt\it $\cal G\mskip-2mu$eometry \&\ 
$\cal T\!\!$opology $\cal P\!$ublications}}  

\def\recd{{\small Received:\qua\receiveddate\ifx\reviseddate\relax
\else\qquad Revised:\qua\reviseddate\fi\par}} 

\def\lognumber#1{\def\thelognumber{#1}}
\def\volumenumber#1{\def\thevolumenumber{#1}}
\def\volumeyear#1{\def\thevolumeyear{#1}}
\def\papernumber#1{\def\thepapernumber{#1}}
\def\pagenumbers#1#2{\def\startpage{#1}\def\finishpage{#2}}
\def\published#1{\def\publishdate{#1}}

\def\received#1{\def\receiveddate{#1}}

\def\accepted#1{\def\accepteddate{#1}}

\long\def\asciiabstract#1{\long\def\theasciiabstract{#1}}

\let\\\par\let\thelognumber\relax\let\thevolumenumber\relax
\let\thepapernumber\relax\let\thevolumeyear\relax\let\startpage\relax
\let\finishpage\relax\let\publishdate\relax\let\receiveddate\relax
\let\reviseddate\relax\let\accepteddate\relax\let\theasciititle\relax
\let\theasciiauthors\relax
\let\theasciiabstract\relax

\let\theasciiemail\relax

\ifplaintex
\font\logobig=cmssbx10 scaled 3836
\font\logomed=cmssbx10 scaled 2557
\else
\font\logobig=cmssbx10 scaled 4200
\font\logomed=cmssbx10 scaled 2800
\fi

\long\def\makeagttitle{   
\count0=\startpage
\agt\hfill      
\hbox to 45truept{\vbox to 0pt{\vglue -13truept{\logomed A\kern -.37em{\logobig 
T}\kern -.38em G}\vss}\hss}
\break
{\small Volume \thevolumenumber\ (\thevolumeyear)
\startpage--\finishpage\nl
Published: \publishdate}

\vglue .25truein

{\parskip=0pt\leftskip 0pt plus
1fil\def\\{\par\smallskip}{\Large\bf\thetitle}\par\medskip} \vglue
0.05truein

{\parskip=0pt\leftskip 0pt plus 1fil\def\\{\par}{\sc\theauthors}
\par\medskip}%
 
\vglue 0.03truein 

{\small\leftskip 25truept\rightskip 25truept{\bf Abstract}\stdspace\theabstract

{\bf AMS Classification}\stdspace\theprimaryclass
\ifx\thesecondaryclass\relax\else; \thesecondaryclass\fi\par
{\bf Keywords}\stdspace \thekeywords\par}\vglue 7truept

}   

\ifplaintex
\font\phead=cmsl9 scaled 950
\font\pnum=cmbx10 scaled 913
\font\pfoot=cmsl9 scaled 950
\headline{\vbox to 0pt{\vskip -4.5mm\line{\small\phead\ifnum
\count0=\startpage ISSN 1472-2739 (on-line) 1472-2747 (printed)
\hfill {\pnum\folio}\else\ifodd\count0\def\\{ }%
\ifx\theshorttitle\relax\thetitle\else\theshorttitle\fi\hfill{\pnum\folio}
\else\def\\{ and }{\pnum\folio}\hfill\ifx\theshortauthors\relax\theauthors
\else\theshortauthors\fi\fi\fi}\vss}}
\footline{\vbox to 0pt{\vglue 0mm\line{\small\pfoot\ifnum\count0=\startpage
\copyright\ \gtp\hfill\else
\agt, Volume \thevolumenumber\ (\thevolumeyear)\hfill\fi}\vss}}
\else
\font=cmsl9 scaled 1050
\font\lnum=cmbx10 
\font=cmsl9 scaled 1050
\makeatletter
\def\@oddhead{{\small\lhead\ifnum\count0=\startpage ISSN 1472-2739 
(on-line) 1472-2747 (printed)\hfill {\lnum\number\count0}\else\ifodd\count0
\def\\{ }\ifx\theshorttitle\relax \thetitle \else\theshorttitle\fi\hfill
{\lnum\number\count0}\else\def\\{ and }{\lnum\number\count0}
\hfill\ifx\theshortauthors\relax 
\theauthors\else\theshortauthors\fi\fi\fi}}\def\@evenhead{\@oddhead}
\def\@oddfoot{\small\lfoot\ifnum\count0=\startpage\copyright\ \gtp\hfill\else
\agt, Volume \thevolumenumber\ (\thevolumeyear)\hfill\fi}
\def\@evenfoot{\@oddfoot}
\makeatother
\fi
\let\maketitlepage\makeagttitle
\let\makeshorttitle\maketitlepage
\let\maketitle\maketitlepage

\newwrite\gtoutfile
\long\gdef\makeheadfile{  
{\def\\{, }\def\s{ }
\immediate\openout\gtoutfile head.xxx
\immediate\write\gtoutfile{To: math@arxiv.org}
\immediate\write\gtoutfile{Subject: put OR rep NNNNN:ppppp}
\immediate\write\gtoutfile{--text follows this line--}
\immediate\write\gtoutfile{Proxy-for: \ifx\theasciiauthors\relax
\theauthors\else\theasciiauthors\fi\s<\ifx\theasciiemail\relax\theemail\else\theasciiemail\fi>}
\immediate\write\gtoutfile{\noexpand\\}
\immediate\write\gtoutfile{Authors: \ifx\theasciiauthors\relax
\theauthors\else\theasciiauthors\fi}
{\def\\{ }\immediate\write\gtoutfile{Title: \ifx\theasciititle\relax
\thetitle\else\theasciititle\fi}}
\immediate\write\gtoutfile{Subj-class: GT or SG, GR etc}
\immediate\write\gtoutfile{MSC-class: \theprimaryclass\ifx\thesecondaryclass\relax\else, \thesecondaryclass\fi}
\immediate\write\gtoutfile{Journal-ref: Algebr. Geom. Topol. \thevolumenumber\s
(\thevolumeyear) \startpage-\finishpage}
\immediate\write\gtoutfile{Comments: Published by Algebraic and
Geometric Topology at}
\immediate\write\gtoutfile{\s\s\s  http://www.maths.warwick.ac.uk/agt/AGTVol\thevolumenumber/agt-\thevolumenumber-\thepapernumber.abs.html}
\immediate\write\gtoutfile{\noexpand\\}
\immediate\write\gtoutfile{}
\ifx\theasciiabstract\relax
\immediate\write\gtoutfile{\theabstract}\else
\immediate\write\gtoutfile{\theasciiabstract}\fi
\immediate\write\gtoutfile{}
\immediate\write\gtoutfile{\noexpand\\}
\immediate\write\gtoutfile{}
\immediate\closeout\gtoutfile}}  

\def\maketitlepage{\makeagttitle\makeheadfile}
\let\makeshorttitle\maketitlepage
\let\maketitle\maketitlepage

\lognumber{32}
\volumenumber{2}
\volumeyear{2002}
\papernumber{32}
\published{2 October 2002}
\pagenumbers{757}{789}
\received{9 April 2002}
\accepted{9 September 2002}

\usepackage{amsmath,amssymb,amscd,graphicx}
\usepackage[all]{xy}

\newtheorem{lem}{Lemma}[section]
\newtheorem{prop}{Proposition}[section]
\newtheorem{thrm}{Theorem}[section]
\newtheorem{cor}{Corollary}[section]
\newtheorem{tm}{Theorem}
\newtheorem{crrr}{Corollary}

\newcommand{\R}{{\bf{R}}}

\newcommand{\Z}{{\bf{Z}}}

\newcommand{\SSS}{{\Sigma}}
\newcommand{\sig}{{\sigma}}

\newcommand{\hfp}{{{HF^+}}}
\newcommand{\hfm}{{{HF^-}}}

\newcommand{\hfpm}{{{HF^{\pm } }}}
\newcommand{\hfi}{{{HF^{\infty} }}}
\newcommand{\hfhat}{{{\widehat{HF}}}}
\newcommand{\cfp}{{{CF^+}}}
\newcommand{\cfm}{{{CF^-}}}
\newcommand{\cfpm}{{{CF^{\pm } }}}
\newcommand{\cfhat}{{{\widehat{CF}}}}
\newcommand{\cfi}{{{CF^{\infty} }}}
\newcommand{\cfs}{{{\widehat{CF} _s}}}

\newcommand{\hfps}{{{HF_s^+}}}
\newcommand{\cfps}{{{CF_s^+}}}
\newcommand{\hfms}{{{HF_s^-}}}
\newcommand{\cfms}{{{CF_s^-}}}
\newcommand{\Zu}{{{\Z [u^{-1}]}}}
\newcommand{\spi}{{{\mathbf s}}}
\newcommand{\spinc}{{{\text{Spin}^c}}}
\newcommand{\gr}{{{\text{gr} }}}

\newcommand{\spr}{{{ \sigma '}}}
\newcommand{\OS}{{{Ozsv\'ath-Szab\'o}}}
\DeclareMathOperator{\imm}{im}

\begin{document}

\title{Floer homology of surgeries on two-bridge knots}
\author{Jacob Rasmussen}
\address{Department of Mathematics, Harvard University\\Cambridge, 
MA 02138, USA}
\email{jacobr@math.harvard.edu}

\begin {abstract}
 We compute the \OS \ Floer homologies \( \hfpm\)
and \( \hfhat \) for three-manifolds obtained by integer
 surgery on a two-bridge knot.
\end{abstract}

\asciiabstract{
We compute the Ozsvath-Szabo Floer homologies HF^{+-} and HF-hat for
three-manifolds obtained by integer surgery on a two-bridge knot.}

\maketitle

\section{Introduction}

The groups \( \hfpm \)  introduced by Ozsv\'ath and Szab\'o in \cite{OS1},
\cite{OS2} have shed new light on our understanding of Floer homology
for three-manifolds.
Conjectured to be equal to the equivariant Seiberg-Witten Floer groups, 
the Ozsv\'ath-Szab\'o groups have most of the known properties of 
these groups, as 
well as some (such as the exact triangle) which were only conjectured. 
In addition, they appear to be more computable than the Seiberg-Witten Floer
groups; at the least, there is an algorithmic procedure to find the
generators for the chain complex associated to a given three-manifold. 
On the other hand, there are some potential stumbling blocks: the size of
the  chain complex is typically much larger than that of the associated 
homology groups, and the differentials in the complex can be difficult to 
determine. 

In this paper, we compute the Ozsv\'ath-Szab\'o Floer homology for integral
 surgeries on
two-bridge knots. Although these are some of the simplest available 
three-manifolds, the fact that the corresponding Seiberg-Witten Floer
 groups are still unknown indicates the computational effectiveness of 
the Ozsv\'ath-Szab\'o groups. Moreover, the method of calculation provides
 some grounds for optimism: the formal properties of their chain complex
enable us  to compute its homology
without having to understand all of its  differentials. 

The result of our computation may be summarized as follows: all the Floer 
homologies associated to a two-bridge knot \(K\) are determined by two 
classical invariants of the knot --- the Alexander polynomial and the 
signature \( \sig (K) \). For the moment, we restrict ourselves to 
describing some interesting special cases. 

\begin{tm}
Let \( K \) be a two-bridge knot and \(T\) the \((2,2n+1)\) torus knot of
the same signature. Let \( K^0 \) and \(T^0\) be their 0-surgeries. 
Denote by \( \spi _k \) the \( \spinc \) structure 
on \( K^0 \) with \(c_1 (\spi _k) = 2k \). Then  
\( \hfp (K^0, \spi _k ) \cong Q \oplus \hfp (T^0, \spi_k) \), 
where \(Q \) is a 
free \( \Z \) module concentrated in a single grading. 
\end{tm}

\noindent As a corollary, we find that 
the Ozsv\'ath-Szab\'o analogue of Fr{\o}yshov's 
$h$-inv\-ar\-iant is determined by the signature:

\begin{crrr}
Let \( K \) be a two-bridge knot and \( K^1 \) the manifold obtained 
by 1-surgery on it. Then \( d(K^1) =
 \min (0,- 2 \lceil \sigma (K) /4 \rceil )\).
\end{crrr}

We now give a quick overview of the calculation. 
 The first two sections of the paper are
 an elaboration of section 8  of \cite{OS1}, which computes \( \hfpm \) and
\( \hfhat \) for the simplest two-bridge knots --- the \((2,n) \) torus 
knots. First, any surgery on a two-bridge knot \(K\) 
 admits a natural Heegaard splitting of genus 2.
 This is described in section 2. In section 3,  we consider \(n\)
 surgery on \(K \) for \( n \gg 0 \). In this case, the generators of
 the chain complex \( \cfhat (K^n) \) are particularly easy to write down.
 Using
certain  ``annular differentials'' described in \cite {OS1}, we compute
the grading in this complex. In section 4, we use these annular differentials
to compute \( \hfhat (K^n) \) for a particular \( \spinc \) structure.
 We next consider the  complex \( \cfp \). Although we do not know 
what most of the differentials in this complex are,
its symmetries allow us to compute \( \hfp \) in every 
\( \spinc \) structure, using only the one group we computed directly.
In section 5 we  apply the exact triangle of \cite{OS2} to get 
\( \hfp (K^n ) \) for any \(n \in \Z \). Finally, in section 6 we 
briefly discuss the extent to which the methods of this paper apply to 
other knots. 

The author would like to thank Peter Kronheimer for his advice and support, 
Kim Fr{\o}yshov for helpful discussions regarding the {\it h}-invariant,
and Peter Ozsv\'ath and Zolt\'an Szab\'o for their interest in the
subject.

\section{Two-bridge knots}

We begin by reviewing some classical facts about two-bridge knots. Good 
references for this material are \cite{Mur} and \cite{BZ}. First, the 
definition:
a  two-bridge knot is any knot \( K\) which admits a presentation of
 the 
form shown in Figure 1: two maxima, two minima, and a braid \(B_K \)
in between. 
If we split \( S^3 \) along a horizontal plane, we can write the 
pair \( (S^3,K) \) as the union of two copies of \( (B^3, \ell _1 \cup 
\ell _2 )\), where \(\ell _1 \) and \(  \ell _2 \) are a pair of parallel
line segments. The two pieces are  glued together by the element of 
the mapping class group of \(S^2 \) with four marked points
 which corresponds to the braid \(B_K \). 

\begin{figure}[ht!]
\centerline{\includegraphics{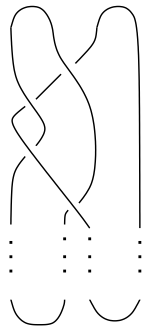}}
\caption{A two-bridge knot}
\end{figure}

Since the double cover of \(B^3 \) branched along \( \ell _1 \cup \ell _2 \)
is a solid torus, the double cover of \( S^3 \) branched along
\(K \) is a lens space. The basic fact about two bridge knots is that
they are classified by their branched double covers.

\begin{prop}[Schubert] 
For every oriented lens space \( L(p,q) \) with \(p\) odd,
 there is a unique two-bridge
knot \(K(p,q) \) with branched double cover \(L(p,q) \). 
\end{prop} 

{\bf Remarks}

\begin{enumerate}

\item The spaces \(L(p,q) \) with \(p\) even  are obtained as branched 
double covers of two-bridge {\it links}. 

\item Lens spaces which are orientation preserving diffeomorphic have the
same knot, so that {\it e.g.} \( K(p,q) = K(p,q') \) if \( q' \equiv  q^{-1}
\mod p \). Lens spaces which are orientation {\it reversing} diffeomorphic
correspond to mirror image knots. In particular, \( K(p,q) \) is amphichiral
if and only if \( q ^2 \equiv -1 \mod p \). 

\item The braid \( B_{K(p,q)} \) can be obtained from a continued 
fraction expansion of \(p/q \), as explained in \cite{Mur}.

\end{enumerate}

Some simple examples of two-bridge knots include the left-hand
\( (2,p) \) torus knots, which are \( K(p,1) \) in our notation, and the 
twist knots \( K(p,2) = K(p,\)\break\((p+1)/2 ) \). 

\subsection{Schubert normal form}
 The knot \(K(p,q) \) admits a canonical projection known as the 
Schubert normal form,  which is 
particularly well-adapted to our purposes. In this section, we describe this 
form and its relevant properties. Proofs may be found in 
\cite{BZ} or \cite{Shu}.

\begin{figure}[ht!]
\centerline{\includegraphics[width=3.6in]{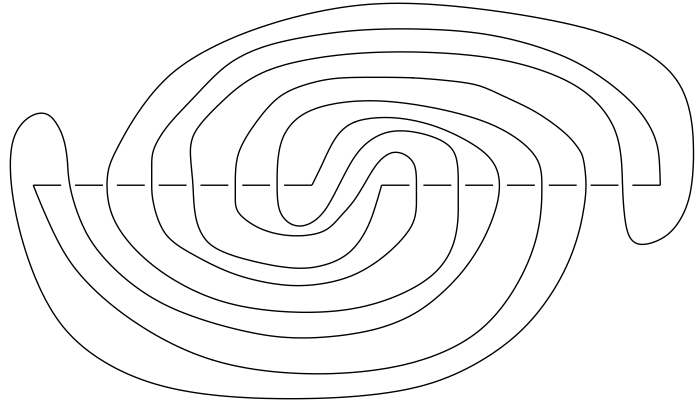}}
\caption{\label{fig:schuex}
The Schubert normal form of \(K(7,3) \)}
\end{figure}

 Consider the knot \(K(p,q) \). Without loss of 
generality, we can assume \(-p < q < p \) and that \(q\) is odd.
The normal form of  \( K(p,q) \) is a projection 
of the sort shown in Figure~\ref{fig:schuex}.
 We break \(K(p,q) \) into 4 segments:
two {\it underbridges} \(U_1 \) and \(U_2 \) (drawn horizontally) 
and two { \it overbridges}
\(O_1 \) and \(O_2 \). If we travel along one of the
 \(U_i \), we go through \( p -1 \) undercrossings, which alternate between
\(O_1\) and \(O_2 \). Similarly, if we  travel along one of the \(O_i \), we
go through \(p-1 \) overcrossings, alternating between \(U_1 \)
and \(U_2 \). 

\begin{figure}[ht!]
\centerline{\includegraphics[width=3.2in]{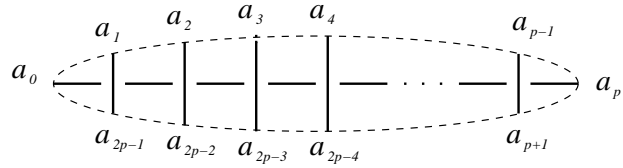}}
\caption{\label{fig:underbridges}
A neighborhood of an underbridge}
\end{figure}

To give a precise description, we consider small neighborhoods of the
underbridges. One such neighborhood is shown in 
Figure~\ref{fig:underbridges}. Abstractly, it is a disk with 
marked points \( a_0, a_1, \ldots a_{2p-1} \) on its boundary. To get the 
Schubert normal form, we glue these two disks together by an orientation
reversing diffeomorphism of \(S^1 \) which identifies the point \(a_{i}
\) on one disk with the point \( a_{q-i} \) on the other. (All the labeling
is modulo \(2p \).) The resulting diagram is most naturally thought of
as living on a sphere, but if we want to draw it, we project onto
a plane.

The Schubert normal form is not quite unique, since
\( K(p,q) = K(p, q^{-1}) \). For most knots (those with \(q^2 \not \equiv
\pm 1 \mod p \)), there are two potential choices of normal form, which are 
related as follows. Suppose we start with the diagram for \(K(p,q) \), 
which we think of as living in a plane. We can straighten out the 
overbridges by a plane isotopy, but at the cost of twisting up the 
underbridges. 
Flipping the resulting diagram over gives the normal form for
\(K(p,q^{-1}) \). 

\subsection{ A Heegard splitting for \(S^3 - K(p,q)\)}

The complement of a two-bridge knot admits a nice handle decomposition:

\begin{prop}
\label{handle1}
 \(S^3 - K(p,q) \) is the union of a genus two handlebody and
a two-handle attached along a curve \( \beta _{p,q} \). This curve is the 
boundary of a regular neighborhood of an overbridge in the Schubert normal
form of \(K(p,q) \).  
\end{prop}

\begin{proof}
 To obtain the desired decomposition,
 we start with the Schubert normal form of \(K(p,q)\). The overbridges 
lie in a plane and the underbridges dip down below them. Split
 \(S^3 \) along this plane, and remove
tubular neighborhoods of  the underbridges. The  part of
 \(S^3 \)   which remains below the plane  is a  
handlebody of genus 2. To
get the part of \(S^3 -K(p,q) \) above the plane, we  must separate
the overbridges by adding a two-handle along   \( \beta _{p,q} \).
Note that \(\beta _{p,q} \) does not depend on which overbridge we choose,
since the complement of a regular neighborhood of one overbridge is a 
regular neighborhood of the other.
\end{proof}

A sample handle decomposition is shown in Figure 4,
following the conventions of \cite{OS1},
Section 8. The region shown is part of the genus 2 surface \( \SSS _2 \)
 which bounds the  handlebody. The remainder of \( \SSS _2 \) consists
of two tubes below the plane of the picture, joining \(A_1 \) to
\(A_2 \) and \(B_1 \) to \(B_2 \). The attaching circles of the handlebody
 are the two horizontal lines \( \alpha _1 \) and \( \alpha _2 \).
The curve \( \beta _{p,q} \) is shown in bold. It intersects each \(
\alpha _i \) \(p \) times. 

\begin{figure}[ht!]
\centerline{\includegraphics[width=\hsize]{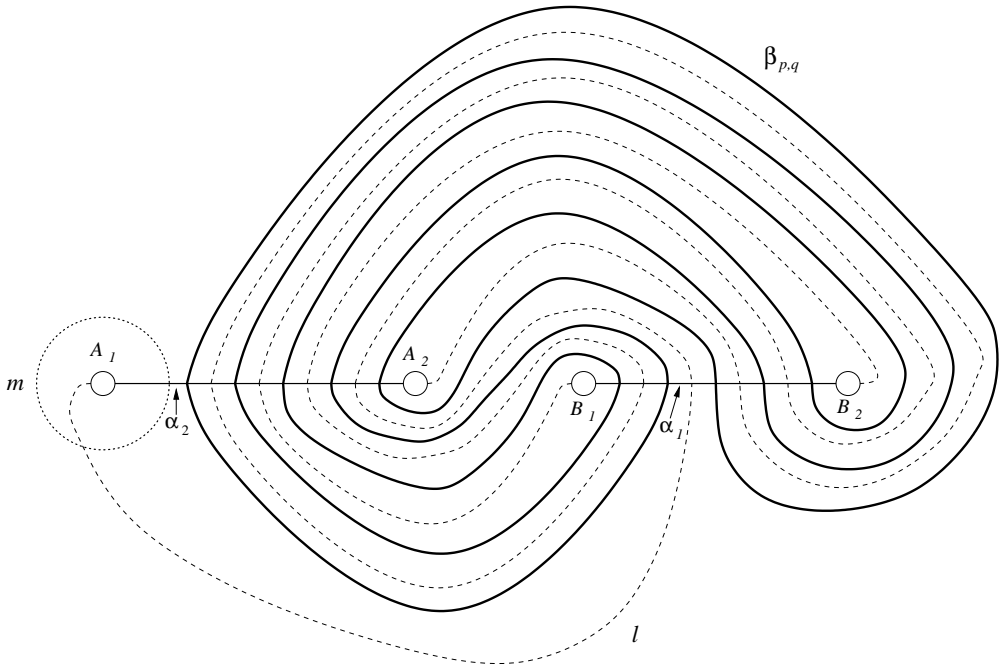}}
\caption{\label{fig:complement}
The Heegaard splitting for the complement of \(K(5,3) \) (the 
figure-eight knot.) The longitude and meridian are shown by dashed lines.
In a general knot, the longitude will spiral
 around the hole labeled 
\(A_1 \), so that it has intersection number
zero with the class of \( \alpha _1 - \alpha _2 \).   }
\end{figure}

\section{The generators of \(\cfhat (K^n) \)}

From now on, we suppose that \(K \) is a two-bridge knot and drop the 
\((p,q)\) when they are not relevant. We want to compute the various
Ozsvath-Szabo Floer homologies for the three-manifold obtained by
\(n\) surgery on \(K\), which we denote by \(K^n\). Until the last
section, we will consider only the case \(n \gg 0 \).

Our aim in this section is to write down generators for the complex 
\( \cfhat (K^n) \) and describe their relative gradings. Of course, 
these gradings
will vary with the choice of \( \spinc \) structure on \( K^n \). 
There is a natural
way to label these \( \spinc \) structures as \( \spi _k \), with respect to 
which we have the following result.

\begin{prop}
\label{prop:grading}
For \( n \) sufficiently large, \( \cfhat (K^n(p,q), \spi _k ) \)
(\(q\) odd)  admits a 
presentation with \( p \) generators \( x_1, x_2 , \ldots x_p \). This
 presentation is independent of the value of  \(n\). For \( k \geq g(K)
\), the complexes \(\cfhat(K^n, \spi _k )\) are all identical to a single
complex which we refer to as \( \cfs (K) \). The grading in
this complex is given by 
\begin{equation*} 
\gr (x_{i+1}) - \gr (x_{i}) = (-1)^{\lfloor \frac{iq}{p} \rfloor }
\end{equation*} 
and there is a unique differential
between \(x_{2i} \) and \(x_{2i + 1 } \). 
In general, the complex \( \cfhat(K^n, \spi _k) \) is obtained by reflecting
the complex  \( \cfs (K) \) at level \(k\). 
\end{prop}

\noindent The process of {\it reflection} will be described in section 3.3.

\subsection{The Heegaard splitting}
Our  first order of business is to construct a Heegaard splitting for 
\(K^n \). We start with the handlebody decomposition of
 \(S^3 -K(p,q^{-1}) \)
described in Proposition~\ref{handle1}. 
(The reason for our choice of \(q^{-1} \)
rather than \(q \) will be apparent in  3.2.)
 When we do surgery, we attach a two-handle along
a curve in the knot complement, and then fill in with a three-ball. 
Alternately, we can just think of attaching the two-handle to \( \SSS _2
 - \beta _{p,q}   \). From this point of view, 
the homology classes of the longitude and meridian are represented by the
curves \( \ell \) and \( m \) shown in Figure ~\ref{fig:complement}.

\begin{figure}[ht!]
\centerline{\includegraphics[width=\hsize]{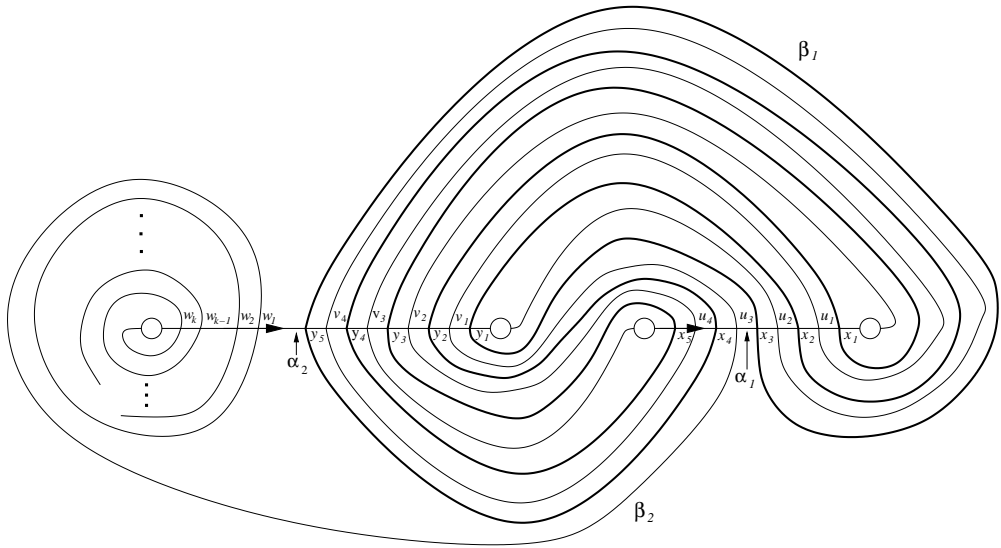}}
\caption{\label{fig:hsplit}
The Heegaard splitting for large \(n\) surgery on \(K(5,3) \).
Note that for for large \(n\), \( \beta _2 \) always has many clockwise
spirals, although the exact number depends on the knot.}
\end{figure}

To do \(n\)-surgery, we attach the two-handle along the curve \(\beta _2
= \ell + n m \). The resulting three-manifold has a Heegaard diagram with 
attaching circles  \( \alpha _1, \alpha _2, \beta _2\), and \( \beta _1 =
\beta _{p,q} \). An example is shown in  Figure~\ref{fig:hsplit}. There
are several general features worth noting. First, the curve \( \beta _1 \)
separates the four-times punctured sphere into two components: the 
{\it exterior}, which contains the point at infinity, and the {\it interior}.
We orient the \( \alpha _i \) so that they both point out of the exterior 
region and into the interior. Thus \( \beta _1 \) has the same intersection
number with each of the \( \alpha _i \). For \( n \gg 0 \), the curve
\( \beta _2 \) always has a region with a large number of clockwise 
spirals, which we draw in the exterior.

Recall that \( \cfhat \) is generated by the intersection points of the 
totally real tori \(T_{\alpha} \), \(T_ {\beta } \) in \( s^2 \SSS _2 \),
or, equivalently, by unordered pairs of intersection points \(\{ x, y \} \)
between the \( \alpha _i \) and the \( \beta _i \). We distinguish
five kinds of such intersection points: between \( \alpha _1 \) and 
\(\beta _1 \), \( \alpha _2 \) and \( \beta _1 \), \( \alpha _1 \) and
\( \beta _2 \), \( \alpha _2 \) and \( \beta _2 \) (not in the spiral),
and \( \alpha _2 \) and \( \beta _2 \) (in the spiral). We label these
points by \(x_i, y_j, u_k, v_l, w_m \) respectively, as shown in 
Figure~\ref{fig:hsplit}. (The convention is that the segment of 
\( \alpha _1 \) between \(x_1 \) and the hole lies in the 
interior region.)  {\it A priori,} we are concerned with pairs of 
intersection points of the form \( \{x_i, v_l \}, \{x_i, w_m \} \)
and \( \{ y_j, u_k \} \). 

\subsection{ \( \epsilon \)--grading and basepoints }

Recall from \cite{OS1} that there is an affine grading
\begin{equation*}
 \epsilon \co T _\alpha \cap T _\beta  \to \text{\it Affine} (H_1 (K^n)) .
\end{equation*}
 If we fix
a basepoint \( z \in \SSS _2 \), this grading corresponds to the usual 
grading of Seiberg-Witten Floer homology by \( \spinc \) structures.
Often, however, it is more convenient to fix an \( \epsilon \)--equivalence
class and vary the basepoint to get all 
\( \spinc \) structures. This is the approach we will take. 

Our first observation is that we need only consider the pairs
\(\{ x_i, w_m \} \):

\begin{lem}
\label{lem:noys}
For \(n \) sufficiently large, there is an \( \epsilon \)--equivalence
 class which contains only points of the form
 \( \{x_i, w_m \} \).
\end{lem}

\begin{proof} There are \(n \) equivalence classes, but the number of
pairs \(\{x_i, v_l \} \) and \( \{ y _j, u _k \} \) is bounded independent
of \(n \).
\end{proof}

To describe \( \epsilon \) explicitly, we need to discuss \( H_1 (K^n ) \).
We start with \(H_1 (\Sigma _2 ) \), which is a free abelian group generated
by elements \( A_1, A_2, B_1, B_2 \), where \( A_i \) is the class
represented by the curve \( \alpha _i \), and the \( B_i \) link the 
``holes'' in Figure~\ref{fig:hsplit},
 so that \( A_i \cdot B_j = \delta _{ij} \). 
We have
\begin{align*}
H_1 (K^n)  & \cong H_1 (\Sigma _2 ) / \langle \alpha _1, \alpha _2, \beta _1,
 \beta _2 \rangle \\
 & \cong \Z  \langle B_1 \rangle \oplus \Z \langle B_2 \rangle
 / \langle B_1 + B_2, \beta _2 \rangle \\
 & \cong \Z \langle B_2 \rangle  / \langle n B_2 \rangle.
\end{align*}
It  follows that the quotient map
\( H_1 (\SSS _2 ) \to H_1 (K^n ) \) is given by\nl
\hbox{}\hspace{1.6in}\( x \to [x \cdot (A_1 - A_2 )] B_2 \).

\begin{lem} There are affine gradings
\begin{align*}
 \epsilon _x \co  &  \{x_i\}  \to \text{\it Affine}( \Z) \\
 \epsilon _w \co &  \{w_m\} \to \text{\it Affine}(\Z / n) 
\end{align*}
so that \( \epsilon(\{ x_i,w_l\}) - \epsilon(\{x_j,w_m \}) =
 ([\epsilon _{x}(x_i) - \epsilon _{x}(x_j) ] + [\epsilon _{w}(w_l)
  - \epsilon _{w} (w_m) ]) B_2 \).
\end{lem}

\begin{proof} (Following Lemma 8.3 of \cite{OS1}.) To define
\( \epsilon_{x} (x_i) - \epsilon_{x} (x_j) \), join \(x_i \) to \(x_j \)
by a path \(a \) along \( \alpha _ 1 \) and a path \(b\) along 
\( \beta _1 \). Then for any \(l \), \( \epsilon \{x_i, w_l \} -
\epsilon \{x_j, w_l \} \) is the homology class of the path \(a-b \),
{\it i.e.} \( [(a-b) \cdot (A_1 - A_2) ] B_2 \). 
We set
\begin{equation*}
 \epsilon_{x} (x_i) - \epsilon_{x} (x_j) =
 (a-b) \cdot (A_1 - A_2).
\end{equation*}
Note that the right-hand side is independent of our choice of \(a\) and 
\(b\), since
\begin{equation*}
 \alpha _1 \cdot ( A_1 - A_2 ) =  \beta _1 \cdot ( A_1 - A_2 ) = 0 .
\end{equation*}
 The definition of
\(\epsilon _w\) is analogous, but  
\(\beta _2 \cdot (A_1 - A_2) = n \), so \( \epsilon_{w}\)
is only defined modulo \(n\). The lemma now follows from
 the additivity of \( \epsilon \). 
\end{proof}

\begin{lem}
\label{lem:ew}
 \( \epsilon _{w} (w_{m+1}) - \epsilon _{w} (w_{m})  \equiv 1
\mod n\).
\end{lem}

\begin{proof}
This is obvious from Figure~\ref{fig:hsplit}. 
\end{proof}

\begin{lem}
\label{ecalc}
 \( \epsilon _x (x_{i+1}) - \epsilon _x (x_{i}) 
= (-1)^{\lfloor \frac{iq}{p} \rfloor }. \) 
\end{lem}

\begin{figure}[ht!]
\centerline{\includegraphics[width=0.9\hsize]{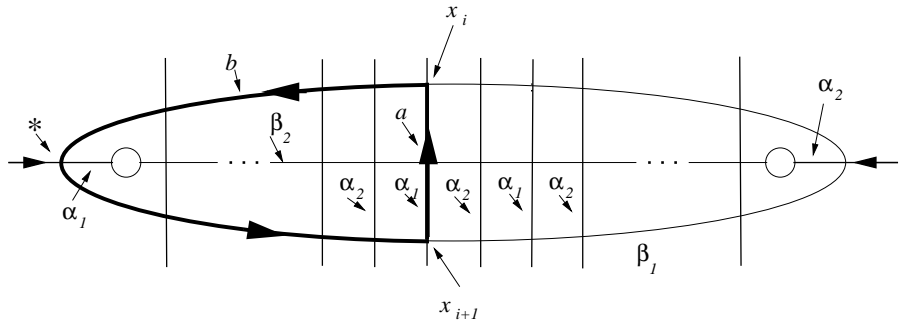}}
\caption{\label{fig:over}
The interior region. We have drawn the case when \(x_i \) is 
above \(x_{i+1}\), but the reverse may also be true. }
\end{figure}

\begin{proof}
Suppose that \(i \) is even, so that the segment of \(\alpha _1 \) between 
\(x_i \) and \( x_{i+1} \) lies in the interior region, which is essentially
a regular neighborhood of an overbridge. 
If we straighten this region out, we get a diagram that looks like
Figure~\ref{fig:over}. 
We choose paths \( a \) and \(b \) as shown. All the contributions
to the intersection number \( (a-b) \cdot (A_1 - A_2 ) \) cancel in pairs,
with the exception of the single point marked with a star. Thus 
\( \epsilon _x (x_{i+1}) - \epsilon _x (x_{i}) \) is \(+1\) if \({x_i} \)
is above \( x_{i+1} \) in the diagram, and \(-1\) otherwise. To see
which case holds, recall from section 2 that straightening out the 
overbridges in the normal form for \(K(p, q^{-1} ) \) gives us the normal
form for \( K(p,q) \), with \( \alpha _1 \) playing the role of an 
overbridge \(O_1 \). The segment of \(\alpha _1 \) between
\(x_i \) and \(x_{i+1} \) corresponds to \(O_1\)'s \(i\)th intersection with 
an underbridge. If we label the intersection points on \( \beta _1 \)
as in Figure~\ref{fig:underbridges}, 
\(x_i \) will correspond to the point \( a_{iq} \), which is in the top 
half of the diagram if and only if \( 0 < iq < p \mod 2p \). 
This proves the result when \(i \) is even. When \(i\) is odd, we 
argue similarly using the exterior region. 
\end{proof}

\noindent Using the formula of Lemma~\ref{ecalc},
 it is not difficult to see that 
\begin{equation*}
 \epsilon _x(x_{i+1}) - \epsilon _x (x_i) = \epsilon _x (x_{p+1-i}) -
\epsilon _x (x_{p-i}) .
\end{equation*}
It follows that 
\( \epsilon _x \) is  symmetric under the 
involution \( \iota \) which sends \(x_i \mapsto x_{p+1-i} \), in the sense
 that 
\begin{equation*}
 \epsilon _x (x_i) - \epsilon _x (x_{j})  = -
 [ \epsilon _x (\iota (x_{i}))- \epsilon _x ( \iota (x_j ))].
\end{equation*}
 As a result,
there is a natural lift of \( \epsilon _x \) to a \( \Z \) valued map
(rather than just \( \text{\it Affine}( \Z ) \) valued), namely the one which has
\( \epsilon _x (x_i) = - \epsilon _x (\iota(x_i)) \).
We call this lift the {\it Alexander grading} on the \( x_i\).
 (The name is explained by Proposition~\ref{prop:alex}.)

Putting Lemmas~\ref{lem:noys} to~\ref{lem:ew} together, we
 get the following corollary, which is an exact 
analogue of Proposition 8.3 of \cite{OS1}.

\begin{cor}
For \(n \) sufficiently large and 
an appropriate choice of \(M \), there is an 
 \( \epsilon \)--equivalence class containing only the
 pairs \( \{x_i, w_{M-\epsilon _x (x_i)}  \} \). 
\end{cor}

From now on, we choose one such equivalence class and work with it. In
this case, we can drop the \(w\)'s and  refer to the generators 
as \(x_i \) without ambiguity.
To get different \( \spinc \) structures, we vary the basepoint. 
As in \cite{OS1}, we consider choices of basepoint which lie inside the 
spiral region.
 Denote by \( z_i \) a point lying in the rectangular region with
corners at \( w_{i+1} \), \(w_{i} \), \(w_{i-1} \), and
\( w_i \) again. (See Figure~\ref{fig:spiral}). We will call the \( \spinc
\) structure determined by the basepoint \(z_{M-k} \) and 
our fixed \(\epsilon\)--equivalence class \( \spi _k \). 
It follows from Lemma 2.12 of \cite{OS1} that
\(\spi _k \) is independent of which \( \epsilon \)--equivalence class we 
took. Although the \(\spi _k \) may not represent every
 \( \spinc \) structure on \(K^n \), we  we will see that 	
they  include all of the interesting ones. 

\subsection{Annular differentials}

In this section, we describe certain differentials in the chain complex
\(CF ^ + \). As a corollary, we obtain the grading in \( \cfhat(K^n,\spi _k) 
\) for every \( \spi _k \).  The key tool  is the following proposition.

\begin{figure}[ht!]
\centerline{\includegraphics{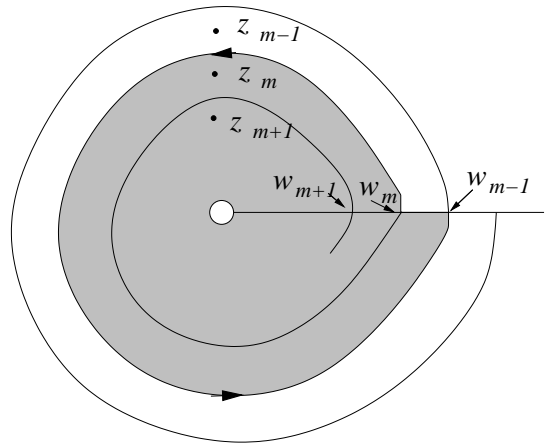}}
\caption{\label{fig:spiral}
A neighborhood of the spiral. The shaded region is part of the 
domain of an interior differential from \( \{x_i,w_{m-1} \} \) to 
 \( \{x_j,w_{m} \} \).} 
\end{figure}

\begin{prop}[Lemma 8.4 of \cite{OS1}]\label{anndif}
 Suppose we are given  \( {\bf x}, {\bf y} \in T _\alpha
\cap T_\beta \subset s^2 \SSS _2 \), and an element
\( \phi \in  \pi _2 ({\bf x},{\bf y}) \) whose domain is an 
annulus with one \( 270^\circ \) corner. Then 
\( \mu ({\bf x}, {\bf y}) = 1 \) and the class \( \phi \) 
has a unique holomorphic representative.
\end{prop}
We refer to such a class as an annular differential from \({\bf x} \) to
\( {\bf y} \), or, if we wish to allow either 
\(\phi \in \pi_2 ({\bf x},{\bf y} ) \) or \(\phi \in \pi_2 ({\bf y},{\bf x}
 ) \), as an annular differential {\it  between} \( {\bf x } \) and
\( {\bf y } \). This definition is useful to us because there is an obvious
 annular differential \( \phi _i \) between \( x_i \) and \(x_{i+1}
\). Its domain can roughly be described as the region bounded by the curves
\(a\) and \(b\) in Figure~\ref{fig:over}.
More precisely, if \( i\) is even, the segment of \(\alpha
_1 \) joining \( x_i \) to \(x_{i+1} \) divides the interior region into
 2 components.  The domain of \( \phi _i \) consists of one of these
components, together with a portion of the spiral like that  shown in 
Figure~\ref{fig:spiral} and the tube joining them. We call \( \phi _i \)
an interior differential. Similarly, if \(i \) is odd, the segment of
 \(\alpha_1 \) joining \( x_i \) to \(x_{i+1} \) divides the exterior
 region into 2 components. The domain of the exterior differential
\(\phi _i \) is the one which contains the spiral, but with the 
shaded region of Figure~\ref{fig:spiral} removed. 

To describe how the annular differentials fit into the complex \( \cfhat \),
we need to fix a \( \spinc \) structure. Initially,
we choose one whose associated basepoint is in the outer part of the spiral.
 More specifically, let 
\(G = \max _i \{ \epsilon _x (x_i ) \} \). (We will show in 3.4 that 
\( G = g(K) \).) Then if \( k \geq G \), the point \( z_{M-k} \) is 
further out in the spiral than all of the \( w_{M-\epsilon _x (x_i )} \). 

\begin{prop}
\label{prop:stablegrading}
  If \( k \geq G \), \( \gr (x_i ) = \epsilon _x (x_i) \) in 
\( \cfhat (K^n, \spi _k ) \). There is a unique differential between
 \( x_{2j} \)  and \( x_{2j+1} \). 
\end{prop}

\begin{proof}
To check the statement about the grading, it suffices to show that 
\begin{equation*}
 \gr(x_{i+1}) - \gr (x_{i}) = \epsilon _x (x_{i+1}) - \epsilon _x (x_{i}) .
\end{equation*}
Suppose \( i \) is even, so that we have an interior  differential
\( \phi _i \) between \( x_i \) and \( x_{i+1} \). Now since
\( k \geq G \), we have \( n_{z_{M-k}} (\phi _i ) = 0 \), so 
\( \gr (x_{i+1}) - \gr (x_{i} ) = \pm \mu (\phi _i) = \pm 1 \), depending
on whether \( \phi _ i \) is a differential from \( x_{i+1} \) to \( x_{i} \)
or the other way around. To tell, we give the boundary of the domain of \( 
\phi _i \) its standard orientation and see which way the \( \alpha \)
segments point. Looking at Figure~\ref{fig:spiral}, we see that the
 differential goes 
from the point with larger \( \epsilon _x \) to the one with smaller
\( \epsilon _x \). This proves the claim.

For \( i \) odd, the argument is similar, except that we now have
\( n_{z_{M-k}} (\phi _i ) = 1 \), and the boundary circle in
 Figure~\ref{fig:spiral} is 
oriented clockwise. Thus the differential goes to the point with larger
\( \epsilon _x \). Without loss of generality, we assume
 \( \epsilon _x(x_{i+1}) <  \epsilon _x(x_{i})\) and compute
\begin{align*}
\gr (x_{i+1} ) - \gr (x_{i} ) & = \mu (\phi _i) -2 n_{z_{M-k}} (\phi _i) \\
& = 1 - 2  \\
& = \epsilon _x ( x_{i+1} ) - \epsilon _x (x_{i}).
\end{align*}
To get the last statement, note that \(n_{x_{M-k}}(\phi _{2j}) = 0 \), 
so the interior differential is ``turned on'' in 
\( \cfhat (K^n, \spi _k ) \). By Proposition~\ref{anndif}, there is a unique
element in  \( {\mathcal M}(\phi _{2j}) \). 
\end{proof}

We now know that the complexes \( \cfhat (K^n, \spi _k ) \) have the 
same generators and grading for \( k \geq G \). We wish to show that 
they have the same differentials as well. To do so, we use the following
result.

\begin{lem}
\label{lem:nindep}
Assume \( n \gg G \), and
suppose \( \phi \in \pi _2 (x_i, x_j) \) has \( \mu (\phi ) = 1 \).
Then one boundary component of \( {\mathcal D }(\phi ) \) consists of 
the segments of \( \alpha _1 \) and \(\beta _2 \) which join \( w_{M - 
\epsilon _x (x_i)} \) and \( w_{M- \epsilon _x (x_j)} \) and which lie
inside the spiral.
\end{lem}

\begin{proof}
Since \( \mu (\phi ) = 1 \), we have 
\begin{equation*}
 n_{z_{M-k}} (\phi ) =  [\gr (x_i) - \gr(x_j) - 1]/2 
\end{equation*}
where the grading is taken with respect to the complex \( \hfhat (K^n, \spi
_k ) \). Since the gradings are constant for \( k \geq G\),
\(n_{z_{M-k}} (\phi ) \) is constant as well. But this can only happen if
\( \partial {\mathcal D }(\phi ) \) is  disjoint from the outer part
 of the 
spiral, and thus from the inner part too. Given this, it is easy to
see that  \( \partial {\mathcal D }(\phi ) \) must have a component as
described. 
\end{proof}

\begin{figure}[ht!]
\centerline{\includegraphics[width=2.5in]{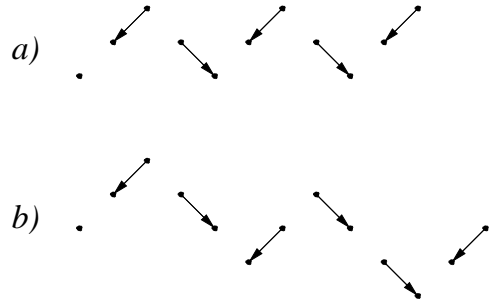}}
\caption{\label{fig:stable}
The stable complexes of (a) \(K(11,5)\) and (b) \(K(13,5)\). The 
generators \( x_1 , x_2, \ldots, x_{p} \) are represented by dots, running
from left to right. The grading of each generator is shown by its height, 
and interior differentials are indicated by arrows.}
\end{figure}

Thus for any differential \( \phi \), 
 we have \(n_{z_k} (\phi ) = n_{z_l} (\phi ) \) if \(k,l \geq G \). 
It follows that for \(n \) sufficiently large and \( k \geq G \), the 
complexes
 \( \cfhat( K^n, \spi _k )\) are all isomorphic to a single complex,
which we call the stable complex of \(K^n \) (written \( \cfs (K^n) \).)
 Although its cohomology is basically trivial, \( \cfs (K^n) \) turns out to
be a very useful object. In particular, if we combine
Lemma~\ref{ecalc} with Proposition~\ref{prop:stablegrading}, we see that the
 grading in \( \cfs (K^n) \) is very easy to compute. 
A few examples of stable complexes  are shown in Figure \ref{fig:stable}.
The reader is encouraged to draw  some others and familiarize himself with
 their behavior.

\medskip 

{\bf Warning}\qua Like 
the Schubert normal form, the stable complex is not quite canonical:
 thinking of \( K (p,q) \) as \( K (p, q^{-1}) \)  gives a different complex.

\medskip 

We now consider what happens to the complex \( \cfhat \) if we choose
a \( \spinc \) structure  \( \spi _k  \) with \( k < G \). We still have
 the same 
interior and exterior differentials \( \phi _i \), but the value of 
\( n_{z_{M-k}}(\phi _i)  \) will change. More precisely, suppose that
 \( \phi \) is 
an interior differential from  \( x_i \) to \( x_j \). Looking
at  Figure~\ref{fig:spiral}, we see that 
we have 
\begin{equation*}
 n_{z_k}( \phi ) = \begin{cases}
 0 & \text{if} \  k< M - \epsilon _x (x_j) \\
1 & \text{if} \  k > M - \epsilon _x (x_i).
\end{cases}
\end{equation*}
Likewise, if \( \phi \) is an exterior differential from 
\( x_i \) to \( x_j \), we have 
\begin{equation*}
 n_{z_k} (\phi ) = \begin{cases}
 1 & \text{if} \ k< M - \epsilon _x (x_j) \\
 0 & \text{if} \ k > M - \epsilon _x (x_i).
\end{cases}
\end{equation*}
Applying this to our differentials \( \phi _i \), we obtain:

\begin{prop}
In \( \cfhat (K^n, \spi _k) \), we have the following proposition.
\begin{equation*}
 \gr (x_{i+1}) - \gr (x_{i} ) = \begin{cases}
 \epsilon _x (x_{i+1}) - \epsilon _x (x_{i}) & \text{if} \ 
\min \{ \epsilon _x (x_i), \epsilon _x (x_{i+1}) \} < k \\
 -[\epsilon _x (x_{i+1}) - \epsilon _x (x_{i})] & \text{if} \ 
\min \{ \epsilon _x (x_i), \epsilon _x (x_{i+1}) \} \geq k.
\end{cases}
\end{equation*}
\end{prop}

\begin{proof}
If \( \min \{ \epsilon _x (x_i), \epsilon _x (x_{i+1}) \} < k \),
then the value of \( n_{z_{M-k}} (\phi _i ) \) is the same as in the 
stable case, so the difference in gradings is the same as well.
 Otherwise, it is the reverse of the stable case
(0 for an exterior differential, 1 for an interior), and we argue 
as in the proof of Proposition~\ref{prop:stablegrading}.
\end{proof}

\begin{figure}[ht!]
\centerline{\includegraphics{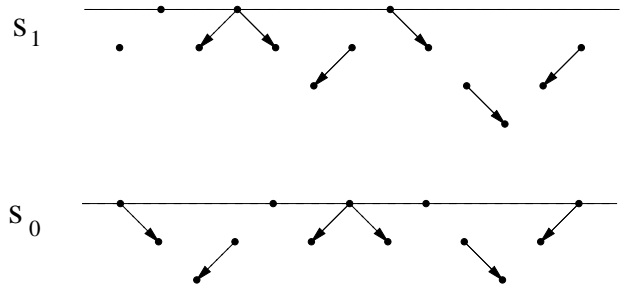}}
\caption{\label{fig:reflect}
 \( \cfhat (K(13,5), \spi _1 ) \) and  
\( \cfhat (K(13,5), \spi _0 ) \). Exterior differentials are shown by 
dashed lines. In each case, the line of reflection is drawn in. (Compare
with Figure \ref{fig:stable}b.) }
\end{figure}

The content of this proposition is  best understood pictorially:
to obtain\break \( \cfhat (K^n, \spi _k ) \), we start with the stable
complex and reflect all the points that lie above level \(k \), as 
illustrated in Figure \ref{fig:reflect}. Note that
in the reflected portion the exterior differentials have been ``turned
on'' and the interior differentials are ``turned off''. 

At this point, we have more or less completed the proof of Proposition
~\ref{prop:grading}. The only thing which remains to be proved is the 
statement that \( \cfhat(K^n, \spi _k ) \) is independent of \(n \) 
for \(n \) sufficiently large. To see this, we use Lemma~\ref{lem:nindep}
again. Indeed, the lemma implies  that the domains of potential
 differentials are independent of 
\(n \). From this, it is not difficult to conclude that the differentials
themselves are independent of \(n\), and thus that the complexes
\( \cfhat (K^n, \spi _k )\) and \( \cfpm (K^n, \spi _k ) \) are as well.
When \(n \gg 0 \), we are thus justified in dropping it from our notation
and simply writing \( \cfhat (K, \spi _k )\) and \( \cfpm (K, \spi _k ) \).

\subsection{Knot invariants and the stable complex}
We end this section by describing some relations between
\( \cfs (K) \) and  classical knot invariants. Unsurprisingly,
there is a connection with the Alexander polynomial:

\begin{prop}
\label{prop:alex}
 Let \(n_k \) denote the number of \(x_i\) in the stable
complex of \( K \) with \( \epsilon _x (x_i) = k \). Then
the  normalized, symmetrized Alexander polynomial of \(K \) is given by 
\begin{equation*}
\Delta _K (t) =  (-1)^{\epsilon _x (x_1)} \sum_k n_k (-t)^k. 
\end{equation*}
\end{prop}

\begin{proof}
From the Heegard splitting of section 2.1, we see that \(\pi _1 (S^3-K) \)
has a presentation \( \langle a_1,a_2 \ | \  w_{\beta_1} \rangle \), where
\(w_{\beta_1} \) is the image of \( \beta_1 \) in \(\pi_1 (H_1) = \langle
a_1,a_2 \rangle \). \( w_{\beta_1} \)
admits the following concrete description: let  \(w \) be the empty word.
Start anywhere on \(\beta _1 \) and travel along it in either direction.
Every time you cross \( \alpha _1 \), append \(a_1^{\pm 1} \) to the end
of \(w\) (\(+1\) if you go up through \(\alpha _1 \) as determined by the 
orientation, \(-1 \) if you go down.) Similarly, append \( a_2^{ \pm 1} \)
to \(w \) every time you pass through \( \alpha _2 \). When you get back 
your starting point, \(w \) will be \( w_ {\beta _1 } \). We 
use the orientation convention of section 3.1, according to which
 the \( \alpha _i \) are oriented so that \([a_1] = -[a_2] = t\) in
 \( H_1 (S^3 - K) \), but any choice of orientations would work just as
well.

We now observe that the Alexander grading is just a geometrical 
interpretation of the free differential, so that   
\begin{equation*}
 \sum (-t) ^{\epsilon _x (x_i)} = \pm t^k  d_{a_1} w_{\beta _1}
\end{equation*}
for some \(k\). Indeed, each term in 
\( d_{a_1} w_  {\beta _1}\) corresponds to an intersection point of 
\( \alpha _1 \) 
and \(\beta _1 \). If \(u\) and \(v\) are two such points, we claim  
that \( \epsilon _x (u) - \epsilon _x (v) \) is the difference between the 
corresponding exponents of \(t\) in \(d_{a_1} w_{\beta _1}\).
 To see this, recall
that \( \epsilon _x (u) - \epsilon _x (v) \) is the homology class in 
\(H_1 (S^3-K) \) of a loop \( \gamma \) which goes from \(u\) to \(v\) along
\(\beta_1 \), and then from \(v \) to \(u\) along \( \alpha _1 \). The class
of \( \gamma \) can be determined by counting (with sign) its intersections
with \( \alpha _1 \) and subtracting off its intersections with \( \alpha
_2 \). The difference in the exponents of \(t\) in \( d_{a_1}w_{\beta _1} \)
counts exactly the same quantity. To check that the signs agree, note that 
the sign of a term in the free differential is the sign of the 
associated intersection point \(u\), which is \((-1)^{\gr (u)}  = (-1)^
{\epsilon _x (u) }\).

Strictly speaking, the Alexander polynomial is the \(\gcd \) of 
\(d_{a_1} w_{\beta _1} \) and \( d_{a_2} w_{\beta _1} \), but an argument of 
McMullen (Theorem 5.1 of \cite{Mc}) shows that when the presentation of
 \( \pi_1 \)
comes from a Heegard splitting, \(d_{a_1} w_{\beta _1} \) and 
 \( d_{a_2} w_{\beta _1} \) agree up to a unit. 
Thus 
\begin{equation*}
\sum (-t) ^{\epsilon _x (x_i)}
\end{equation*}
represents  \( \Delta _K (t)\). It is symmetric,
 so it must be the symmetrized Alexander polynomial.
Finally, it is easy to check
 that our choice of the sign \((-1)^{\epsilon _x(x_1)} \)
normalizes \( \Delta _K \) to have \( \Delta _K (1) = 1\). 
\end{proof}

\begin{cor}
The genus of \(K \) is the maximum value assumed by \( \epsilon _x \).
\end{cor}

\begin{proof}
All two-bridge knots are alternating. For such knots, 
it is well-known that the genus is the highest
power of \(t \) appearing in the symmetrized Alexander polynomial.
\end{proof}

Interestingly, there is also a relationship with the signature.

\begin{prop}
\label{cor:sig}
\( \sigma (K ) = \epsilon _x (x_{1}) - \epsilon _x (x_p ) \)
\end{prop}

\begin{proof}
This is a restatement of Theorem 9.3.6 of \cite{Mur}. 
\end{proof}

\section{Computing the homology}

We now turn to our main task, which is to compute the homology of \( \cfhat
(K, \spi _k ) \) and \( \cfpm (K, \spi _k ) \). As described in the 
introduction, the Floer homology of the zero surgery
 splits naturally into two parts: 
a group \(V\) which looks like the Floer homology for a torus knot of the 
same signature, and a group \(Q \) which is concentrated in a single grading.
In \( \cfhat(K, \spi _k ) \), there will be a third summand \( \Z [u^{-1}]
\). If \( \sig (K) > 0 \), \(V\) will be absorbed into the \( \Z [u^{-1}] \)
 summand, while if \( \sig (K) < 0 \), it will not.

To describe \(Q\) and \(V\) precisely, we need some notation.
 Let \( a_k \) be the coefficient of
 \( t^k \) in the symmetrized Alexander polynomial of \( K \), and set
\begin{equation*}
u _k = \sum_ {i>k} (i-k)a_i.
\end{equation*} 
\( u_k \) is a familiar quantity: it is just
the Seiberg-Witten invariant associated to
the \(k\)th \( \spinc \) structure on \( K^0 \). We let \( \sigma  \)
 be the signature of the two-bridge knot \( K \), and set
\begin{align*}
\spr &= \sigma / 2 \\
 h_k &= \max( \lceil \frac{|\spr | -k}{2} \rceil , 0) \\
 b_k &= u_k + (-1)^{k - \spr }h_k. 
\end{align*}
\noindent Then we have the following:

\begin{thrm}
For \( k \geq 0 \), there is an isomorphism of \( \Z [u ] \) modules
\begin{align*}
 \hfp (K, \spi _k ) \cong  \hfp (K, \spi _{-k} ) \cong
\begin{cases}
 Q_k  \oplus \Z [u^{-1}] & \text{if} \ \sigma \geq 0 \\
  Q_k \oplus V_k \oplus \Z [u^{-1}] & \text{if} \ \sigma \leq 0 
\end{cases}
\end{align*}
where \( Q_k \cong \Z^{|b_k|} \) and \( V_k \cong \Z [u^{-1} ]/ u^{-h_k }\).
All the elements of \(Q_k \) have grading  \( k-1\), and  \( 1\in V_k \) has 
grading \(k -2h_k\) if \( k \equiv \spr \mod 2 \) and
 \(k +1-2h_k \) otherwise. 
Finally, \( 1 \in \Z [u^{-1} ] \) has grading \( \spr \) if \(
\sigma \leq 0\), or \(\spr - 2h_k \) if \( \sigma \geq 0 \).
\end{thrm}
{\bf Warning}\qua For ease of expression, we have chosen to state the gradings 
in absolute terms, although they are  only meaningful relative
to each other. They do {\it not} correspond to the absolute gradings 
 discussed in \cite{OS3}. 

\begin{cor} 
\label{cor:hfhat}
\begin{align*}
\hfhat (K, \spi _k ) \cong 
\begin{cases}
 \widehat{Q_k} \oplus \Z & \text{if} \ \sigma \geq 0 \\
 \widehat{Q_k} \oplus \widehat{V_k} \oplus \Z & \text{if} \ \sigma \leq 0 
\end{cases}
\end{align*}
where \( \widehat{Q_k} \cong Q_k \otimes H^*(S^1) \) and 
\( \widehat{V_k} \cong H^*(S^{2h_k -1}) \). The grading of 
\( 1\in H^*(S^{2h_k -1}) \) is the same as that of \( 1 \in V_k \), 
and the grading of \(1 \in \Z \) is the same as that of \( 1 \in 
\Z [u^{-1}] \). 
\end{cor}

\begin{cor}
\label{cor:hfm}
\begin{align*}
\hfm (K^n, \spi _k ) \cong 
\begin{cases}
  Q_k \oplus \Z [u] & \text{if} \ \sigma \geq 0 \\
  Q_k \oplus V_k \oplus \Z [u] & \text{if} \ \sigma \leq 0 
\end{cases}
\end{align*}
The gradings are shifted from those in \( \hfp \). More precisely, the 
grading of \( 1 \in \Z [u] \) is 2 less than that of \( 1 \in 
\Z [u^{-1}] \subset \hfp ( K^n, \spi _k ) \), while the gradings of 
\( V_k \) and \(Q_k \) are 1 less than those of their counterparts in 
\( \hfp ( K^n, \spi _k ) \).
\end{cor}

\subsection{\( \cfi (K) \) and  \( \cfpm (K) \) } 
Until now, we have only discussed the complex \( \cfhat (K) \). We will
need to use 
 \( \cfi (K), \cfp (K) \), and \( \cfm (K) \) as well, so we collect some
basic facts about them from \cite{OS1}, \cite{OS2} here. 

We begin with \( \cfi \). Visually, this complex may be obtained by 
stacking copies of \(\cfhat \) together, each two units apart from the 
next. The stack extends infinitely in both directions.
 As observed in Theorem 8.9 of \cite{OS1}, the complex we get is 
independent of our choice of \( \spinc \) structure. 
Depending on the situation, we will use one of two different conventions to
label the generators of \(\cfi \). If we have fixed   a \(\spinc \) structure
\( \spi \), we follow the notation of \cite{OS1}
and consider the generators as pairs \([{\bf x}, j ] \). On the other hand,
if we want to avoid specifying a \( \spinc \) structure, it is often 
convenient to choose a zero level for the grading on \( \cfi \) and 
let \(({\bf x}, j)\)
 denote the generator above \( {\bf x} \) with grading \(j\).

For any choice of \( \spinc \) structure \(\spi \), \(\cfi (M) \)
 has a filtration
by subcomplexes  \( \ldots \subset F_{-1}^{\spi} \subset F_0^{\spi}
 \subset F_1^{\spi} \subset 
\ldots \), where \( F_i^{\spi} = \{ [{\bf x},j] \ | \ j < i \} \).
For each \(i\), we have 
\begin{align*}
 \cfm (M, \spi) \cong & F_i^{\spi} \\
\cfp (M, \spi) \cong & \cfi (M) / F_i^{\spi} \\
\cfhat (M, \spi ) \cong & F_i^{\spi}/ F_{i-1}^{\spi}.
\end{align*}
Frequently, we wish to relate \( \hfpm \) to \( \hfhat \). In one direction,
this is accomplished by the following Gysin sequence:

\begin{lem} There is a long exact sequence:
\label{lem:gysin}
$$\longrightarrow\hfhat _i (M, \spi ) \longrightarrow HF_i ^+ (M, \spi ) 
\buildrel{u}\over\longrightarrow HF_{i-2}^+
(M, \spi ) \longrightarrow \hfhat _{i-1} (M, \spi )$$
\end{lem}

\proof
At the chain level, the map \( u : \cfp (M, \spi ) \to \cfp (M, \spi ) \)
is clearly a surjection with kernel \(\cfhat (M, \spi ) \).  
Thus we have a short exact sequence
$$0 \longrightarrow\cfhat_* (M, \spi ) \longrightarrow CF_*^+ (M, \spi ) 
\buildrel{u}\over\longrightarrow CF_{*-2}^+
(M, \spi ) \longrightarrow0\eqno{\qed}$$

Conversely, to go from \( \hfhat \) to \( \hfp \), we can sometimes use
the following:

\begin{lem}
\label{hattoplus}
There is a spectral sequence with \( E_2 \) term \( \hfhat (M, \spi )
\otimes \Z [u^{-1}] \) which converges to 
\( \hfp (M, \spi ) \).
\end{lem}

\begin{proof}
We have a filtration \(\overline{F}_0^{\spi} \subset
 \overline{F}_1^{\spi} \subset \ldots \) of
\( \cfp (M, \spi ) \), where \( \overline{F}_i^{\spi} \) is the image of
 \(F_i^{\spi} \)
under the quotient map. 
The quotient \( \overline{F}_{i+1}^{\spi} / \overline{F}_{i}^{\spi}
 \cong \cfhat (M, \spi ) \), and the usual
arguments show that the first differential is just the differential in 
\( \cfhat \). 
\end{proof}

\subsection{Conjugation symmetry and the antistable complex}

It is well known that there is a  conjugation symmetry
$$ \hfpm (K, \spi _k ) \cong \hfpm (K, \spi _{-k} ). $$ In our case, 
this symmetry is explicitly realized on the chain level by the map
\( \iota \) described in section 3.2. To see this, we observe that 
there is an exact symmetry between the interior and exterior regions of our 
Heegard diagram. 

\begin{lem}
There is an involution \( j : \SSS _2 \to \SSS _2 \) which preserves 
the \( \alpha _i \) and the \( \beta _i \). \( j \) reverses the interior
and exterior regions, and \( j(x _i) = x _{p+1-i}\). 
\end{lem}

\begin{proof}
We begin by defining \( j\) on the punctured \(S^2\) of our usual diagram.
The circle \( \beta _1 \) divides \( S^2 - 4 D^2 \) into the interior and 
exterior regions, each of which is homeomorphic to the punctured
disk of  Figures~\ref{fig:over} and \ref{fig:underbridges}.
 To obtain \(  S^2 - 4 D^2 \), we identify 
2 copies of this disk along their boundary by an orientation reversing 
diffeomorphism which sends the point \(a_i \) to \(a_{q-i} \). Let
\( j \) be the map which interchanges the two copies.
Since \( i \to q-i \) is an involution, \( j \) respects the 
gluing, and thus defines an orientation preserving involution
of \( S^2 - 4 D^2 \). It is easy to see that 
\( j \)  preserves \( \alpha _1, \alpha _2 \)
and \( \beta _1 \), and it would preserve \( \beta _2 \) if there we no 
spiral. We solve this problem by pushing the spiral onto the adjacent tube, 
and extend \( j \) to \( \SSS _2 \) by an orientation preserving
involution which switches the ends of the tubes and is chosen to preserve 
the spiral. Since \(j \)  reverses the ends of \( \alpha _1 \) in 
\( S^2 - 4 D^2 \), it must send  \( x _i \) to \( x_{p+1-i} \).
\end{proof}

It follows that \( \iota \) induces an isomorphism 
\( \cfpm (K, \spi _k) \cong  \cfpm (K, \spi _{-k}) \). (Note, however, 
that we have made different choices of complex structure, {\it etc.} 
 on the two 
sides of the equation. Since we are mostly interested in the homology,
this will not be an issue.) 
When \( k \ll 0\), the complex \( \cfpm (K,\spi _k ) \) is independent of 
\(k\), and we refer to it as the {\it antistable complex} 
\(CF_a^{\pm} (K) \) of \(K\). It is obvious from the discussion above
 that  \( HF_a^{\pm} (K) \cong \hfpm _s (K) \).

\subsection{The homology of the stable complex}

We need one more result before we can work out
\( \hfp (K^n) \):

\begin{prop}
\label{stableh}
 The homology of \( \cfs (K) \) is a single copy of \( \Z \) with the 
same grading as \( x_1 \). 
\end{prop}

At first glance, this seems like a very modest statement.
 In fact, it follows 
easily from the general results of  \cite{OS2}
 that the homology of the stable complex
must be \( \Z \), so only real content of the proposition is to 
tell us the grading of the generator. Nonetheless, this result is a 
key step in the calculation. In particular, it is the only place at
which the differentials make an appearance. 

As is usual with Floer homology, the proof of Proposition~\ref{stableh}
 should be 
thought of in terms of Morse theory. For \(\hfhat \), the difference in the
value of the ``Morse function'' between two critical points is given by the 
area of the Whitney disk joining them, or equivalently by the area
of the associated domain in \( \SSS _2 \).

 Consider the  generators
\( x_{2i} \) and \( x_{2i+1} \) in \( \cfs(K) \). By
 Proposition~\ref{prop:stablegrading}, we know that there is an interior 
differential, which goes either from \(x_{2i} \) to \(x_{2i+1} \)  or
from \( x_{2i+1} \) to 
\(x_{2i}\). To avoid repeatedly stating both possibilities, we introduce 
alternate names \(\hat{x}_{2i} \) and \(\hat{x}_{2i+1} \)
 for \(x_{2i}
\) and \(x_{2i+1} \), arranged so that the differential is always
 from  \(\hat{x}_{2i} \) to \(\hat{x}_{2i+1} \).

 Suppose we choose the metric
on \( \SSS _2 \) so that the interior region is very thin,and thus that the
drop 
in the value of the ``Morse function'' between \(\hat{x}_{2i} \) and 
\( \hat{x}_{2i+1} \)
is  quite small. Now if we were really
doing  Morse theory, the logical conclusion would be that
 \( \hat{x}_{2i} \)
and \( \hat{x}_{2i+1} \) are  a pair of  critical points created by 
some small perturbation, and we should just cancel them.
 If we could do this to all the pairs, we would 
be left with nothing but \(x_1 \), and the result would be obvious. To prove
proposition~\ref{stableh}, we need to  translate this Morse-theoretic
 intuition over to the formal setup of \( \hfhat \).

The first thing to observe is that canceling critical points is not specific
to Morse functions: we can do it in any chain complex. More precisely,
we have the following lemma.

\begin{lem}
\label{lem:cancel}
Suppose \( (C^*, d) \) is a chain complex with \( C^i \) freely generated
by \( x_1, x_2,\ldots x_{j_i} \)  and \( C^{i-1} \) freely generated by
\(y_1, y_2, \ldots y_{j_{i-1}} \). Denote the \(y_j \)th component of 
\(d(x_k) \) by \(d(x_k, y_j)\). Then if \( d(x_1,y_1) = 1 \), \( (C^*,d) \) 
is chain homotopy equivalent to  a new complex
\( (\overline{C} ^*, \overline{d} ) \). The generators of \(
\overline{C}^* \) are the generators of \(C^* \), but with \(x_1 \)
and \(y_1 \) omitted, and \( \overline{d} \) agrees with \(d\) except on
\( \overline{C}^i \), where we have 
\begin{equation*}
\overline{d} (x_i) = d(x_i) + d(x_i,y_1)d(x_1)
\end{equation*}
\end{lem}

\noindent This is most likely a folk theorem. 
A detailed proof may be found in Lemmas 3.7--3.9 of \cite{Fl}. 

To prove the proposition, we want to apply this lemma repeatedly to 
\( \cfs (K) \). To do so, we need to show that we still have 
\( d(\hat{x}_{2i},\hat{x}_{2i+1}) = \pm 1 \) in the new complex.
 This follows from
the fact that we have a ``Morse function.''

Suppose we have chosen a metric  \(g\) and  a basepoint \(z_k\) on
 \( \SSS _2 \). Then for any \(i\) and \(j\), there is a
 unique class  \(\phi ^i _j \in \pi _2
(\hat{x}_i,\hat{x}_j) \) whose domain misses the basepoint. Since
\begin{equation*}
 {\mathcal D}(\phi ^i_k) = {\mathcal D}(\phi ^i_j) +
 {\mathcal D}(\phi ^j_k), 
\end{equation*}
 we can choose a function \( A : \{\hat{x}_i\}  \to \R \) with 
\( A(\hat{x}_i) - A(\hat{x}_j) \) equal to the signed 
area of \( {\mathcal D}(\phi ^i_j) \).
Holomorphic disks always have positive domains, so
\( d(\hat{x}_i,\hat{x}_j) = 0 \) unless \( A(\hat{x}_i) > A(\hat{x}_j) \). 

\begin{lem}
For an appropriate choice of metric on \( \SSS _2 \) the intervals 
\begin{equation*} I_j =
[A(\hat{x}_{2j+1}), A(\hat{x}_{2j})] 
\end{equation*}
are disjoint. 
\end{lem}

\begin{proof}
Since we are working with the stable complex, our basepoint is in the 
outer part of the spiral.
 We think of  the surface \( \SSS _2 \) as divided into 
components by the \( \alpha _i \) and \( \beta _i \) and vary the 
area of each component. Choose the area of the interior region and the
spiral to be some very small number \( \epsilon \).
Thus the length of the interval \( I_j \) is less than \( \epsilon\).
If we choose \( t_j \in I_j \), it suffices to show that we
 can vary the areas of the
exterior components so none of the \( t_j \) are within \( 2 \epsilon \)
of each other.  

To do this, we consider the exterior region, which is divided into \((p+1)/2
\) components by \(\alpha _1 \). We label these components \( C_i \)
(starting with the one nearest to the spiral), and take the area of \(C_i \)
to be
\( M + \eta _i - \eta _{i-1} \), where \( M \gg \eta _i \gg \epsilon
\) (except for \( \eta _0 = \eta _{(p+1)/2} = 0\).) Now the domain of 
\(\phi ^{2k}_{2k-1} \) is the complement of the domain of 
 an exterior differential, 
so the area of \({\mathcal D}(\phi ^{2k}_{2k-1}) \) is a constant
plus or minus \(\eta _{j_k} \) for some \(j_k \) which is uniquely 
associated to \(k\). 
 Thus if  we vary a single \( \eta _j \) while holding the
other \( \eta _i \) fixed, we vary one of the \( t_k - t_{k+1} \)  by
\( \pm \eta _j \) while 
holding the other such differences fixed. By varying first the 
difference \( t_2 - t_1 \), then the difference \( t_3 - t_2 \), and so
forth, we can make sure all the \(t_i \) differ by more than \( \epsilon \).
(This method also shows that we can assume \( A(x_1 ) \) is not
in any of the \(I_j \).)
\end{proof}

\begin{proof}[Proof of Proposition~\ref{stableh}]
 We start with the complex \( \cfs (K) \), which
we label \((C^*_0, d_0 ) \). By Lemma~\ref{lem:cancel},
 we can cancel
\( \hat{x}_2 \) and \( \hat{x}_3 \) to obtain a new complex
 \((C^* _1, d_1) \). We claim
that we still have \( d_1 (\hat{x}_{2i}, \hat{x}_{2i+1} ) = \pm 1 \)
 for \( i > 1\).
Clearly, it suffices to show that either 
\( d_0(\hat{x}_{2i},\hat{x}_3) = 0 \) or
 \( d_0 (\hat{x}_{2},\hat{x}_{2i+1}) \) = 0. But since
\( d_0(x,y) = 0 \) if \( A(y) > A(x) \), this follows 
 from the fact that the intervals \( I_1 \) and \( I_i \) are disjoint.

Moreover, we claim that we still have \(d_1(x,y) = 0 \) if \( A(y ) > A(x)
\). To see this, suppose the contrary. Then since \( d_0 (x,y) = 0 \), we
must have \( d_0(x,\hat{x}_3) \neq 0 \neq d_0(\hat{x}_2,y) \), so
\( A(x) > A(\hat{x}_3) \). In fact  \( A(x) > A(\hat{x}_2) \),
 since the 
intervals \(I _j \) are all disjoint. Similarly, we have 
\( A(\hat{x}_3) > A(y) \), so
\begin{equation*}
A(x) > A(\hat{x}_2) > A(\hat{x}_3) > A(y)
\end{equation*}
which contradicts our assumption that \(A(y) > A(x) \). 

To finish the proof, we simply repeat this process, cancelling 
\(\hat{x}_4 \) and \( \hat{x}_5 \) to obtain a complex \(C_2^* \), and so on
until we reach \( C_{(p-1)/2}^* \), which has \( x_1 \) as its only 
generator. Since all the \( C_i^* \) are chain homotopy equivalent,
the result follows. 
\end{proof}

\begin{cor}
\label{stableh+}
\( \hfps (K) \cong \Z [u^{-1}] \), where the homology class
 \( 1 \in \Z [u^{-1}] \) has the same grading as \( [x_1,0] \in \cfps (K) \).
\end{cor}

\begin{proof}
Follows immediately from the spectral sequence of Lemma~\ref{hattoplus}.
\end{proof}

\begin{cor}
\( \hfms (K) \cong \Z [u] \), where the homology class
\(1 \) has the same grading as \( [x_1, -1] \in \cfms (K) \).
\end{cor}

\subsection{Computing \( \hfp (K ) \)}

To prove Theorem 1,
we consider \( \cfs (K) \) and \( \cfps (K ) \)
  as subsets of \( \cfi (K ) \). Throughout this section, we describe
elements of \( \cfi (K ) \) using the bracket notation
 with respect to a stable \( \spinc \) structure. With this convention,
 the generators of \( \cfs (K) \) are the points \( [x_i, 0 ] \), and the
generators of \(\cfps (K) \) are the points \( [x_i,j] \) with \(j \geq 0 \).
We choose an absolute grading on \( \cfi(K) \) which agrees with \( \epsilon
_x \) on \( \cfs (K) \). From Lemma~\ref{cor:sig}, we see that the 
grading of \( [x_1,0] \) is \( \spr \). 

\begin{figure}[ht!]
\centerline{\includegraphics{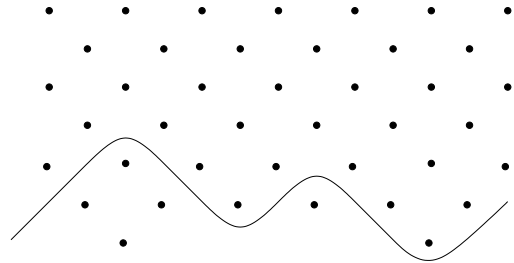}}
\caption{ \label{fig:ck}
The generators of \( \cfp(K(13,5), \spi _0) \). The subcomplex
\(C_0 \) lies below the line, and the quotient complex \( \cfps \) is 
above it.}
\end{figure}

Let us now put \( \cfp (K, \spi_k ) \) into the picture as well. Recall that
to get the generators of \( \cfhat (K, \spi _k) \), we reflect the 
generators of \( \cfs (K ) \) at level \( k \). We view
the former complex as a subset of \( \cfi (K ) \), positioned so
 that the 
unreflected generators (those with \( \epsilon _x \leq k \)) agree with the 
corresponding generators of \(  \cfs (K) \).
 The generators of \( \cfp (K, \spi _k ) \) correspond to 
those generators of \( \cfi (K) \) lying on or above \( \cfhat (K,
 \spi _k ) \). In particular, they include all the generators of 
\( \cfps (K ) \). This situation is illustrated in Figure~\ref{fig:ck}. 

Let \( C_k =  \cfp (K, \spi _k ) \cap \cfms (K) \). Then
 \( C_k \) is a subcomplex of \( \cfp (K, \spi _k ) \), since
\( \cfms (K) \) is a subcomplex of \( \cfi (K) \). Moreover,
if we discard the generators of \( C_k \) from
\( \cfp (K, \spi _k ) \), we are left with \( \cfps (K ) \).
Thus we have a short exact sequence of chain complexes 
\begin{equation*}
\begin{CD}
0 @>>> C_k @>>> \cfp (K, \spi _k ) @>>> \cfps (K) @>>> 0.
\end{CD}
\end{equation*}
We will compute the homology of  \( C_k \) by comparing it to 
\( \cfms (K) \) and \( CF _a^+ (K) \). Since all of the indices in the 
argument can get a bit confusing, we first describe a specific example.
Suppose we wish to compute \( \hfp (K(13,5), \spi _0) \), so the
situation is as shown in Figure~\ref{fig:ck}. \( C_0 \) has generators 
in three rows, with gradings \(0\), \(-1\), and \(-2\). As shown in 
Figure 11a, the upper two rows are identical to the top two 
rows of the complex \( \cfms (K) \). Thus to compute the homology
 in the top row, it suffices to know \( \hfms (K) \). 
Similarly the bottom two rows of \(C_k \) 
are the same as the bottom two rows of \( CF_a^+ (K) \), and we can compute 
the homology of the bottom row of \(C_k \) in terms of \( HF_a^+ (K) \).

We now consider the general situation. 

\begin{lem}
\label{lem:ck}
The set of generators of \( C_k \) is 
\( \{ [x_i, j] \ |\  0 > j \geq k - \epsilon _x (x_i ) \} \).
\end{lem}

\begin{proof}
We get the complex \( \cfhat (K, \spi _k ) \) by reflecting
the complex for \( \cfs (K ) \) at level \( k \). Thus the
generators of \( \cfhat (K, \spi _k ) \subset \cfi (K ) \) are of
the form 
\begin{align*}
& [x_i, k - \epsilon _x (x_i)] 
  & \text{if} \ k \leq \epsilon _x (x_i) \\
& [x_i, 0]  & \text{if} \ k \geq \epsilon _x (x_i)
\end{align*}
Since \( C_k \) is generated by those \( [x_i, j ] \in \cfp (K, \spi _k )
\) with \( j < 0 \), the result follows. 
\end{proof}

\begin{figure}[ht!]
\centerline{\includegraphics[width=3in]{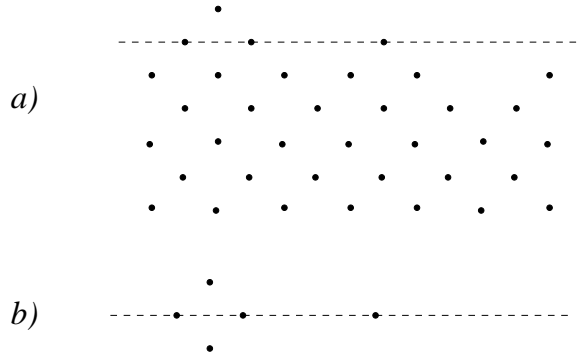}}
\caption{\label{fig:trunc}
To obtain \(C_0 \) for \(K(13,5) \), we (a)
 truncate \( \cfms \) at grading
 \(-1\) and (b) take the union of this set and its reflection.}
\end{figure}

Visually, we can restate this result as follows: to get the complex \( C_k 
\), start with the complex \( \cfms (K ) \), truncate
everything below the horizontal line with grading \( k-1 \), and then add
in the 
reflection about this line, as illustrated in Figure~\ref{fig:trunc}.
Indeed,  truncating at grading
 \( k-1 \) gives all the points of the form \( [x_i, j] \) with \( j < 0 \)
and 
\begin{equation*}
\gr ([x_i,j]) =  \epsilon _x (x_i ) + 2j \geq k-1 
\end{equation*}
or equivalently
\begin{equation*}
 0 > 2j \geq k-1 - \epsilon _x (x_i).
\end{equation*}
Adding in the reflection 
gives us all the points  \( [x_i, j] \) with 
\begin{equation*}
 0 > 2j \geq k-1 - \epsilon _x (x_i) - (-2-[ k-1 -\epsilon _x (x_i)]) = 2
(k- \epsilon _x (x_i ) ).
\end{equation*}
But this is precisely the condition of the lemma.

\begin{cor} 
\label{cor:top}
\( H_i (C_k ) \cong \hfms _{i} (K) \) for
\( i \geq k \).
\end{cor}

\begin{proof} For \( i \) in this range 
 we are in the top half of \( C_k \), where the 
complexes  \( C_k \) and \( \cfms   (K ) \) are identical.
\end{proof}

On the other hand, the reflected half of \(C_k \) is a subcomplex of the 
antistable complex. Indeed, truncating \(\cfs (K) \) above level \(k+1\) and
then reflecting is the same as reflecting \( \cfs (K) \) to obtain 
\( \cfhat _a (K) \) and then truncating below level \(-k-1\). 

\begin{cor}
\label{cor:bottom}
\( H_i (C_k ) \cong \hfps _{i-2k} (K) \) for \( i \leq k-2 \). 
\end{cor}

\begin{proof}
For \(i\) in this range \(C_k\) and \( CF_a^+ (K) \) are the same, so 
\begin{equation*}
 H_i(C_k) \cong HF_{a \ i-2k}^+(K) \cong \hfps _{i-2k} (K) 
\end{equation*}
 by the charge conjugation symmetry of section 4.2. 
\end{proof}

To sum up, we now know \( H_i (C_k ) \) for every \( i \neq k-1 \).
To get \(H_{k-1}(C_k) \) we need an easy result from algebra, whose proof is
left to the reader:

\begin{lem}
\label{lem:free}
Suppose
\begin{equation*}
\begin{CD}
A @>{f}>> B @>{g}>> C \oplus D
\end{CD}
\end{equation*}
is a sequence of free abelian groups with \( gf = 0 \) and 
\( \ker g / \imm f \) free. If \( g' =\pi g : B \to C \) is the obvious
map, then \( \ker g' / \imm f \) is free as well.
\end{lem}

In particular, if we take \( A = \cfms _{k} (K) = (C_k) _{k}\),
\(B = \cfms _{k-1} (K) = (C_k )_{k-1}\), and
\( C \oplus D = (C_k)_{k-2} \oplus D = \cfms _{k-2} (K)\), we see that
\( H_{k-1}(C_k) \) is a free group. 
To compute its rank, we use the  Euler characteristic. 

\begin{lem}
\label{lem:chi}
\begin{equation*}
\chi (C_k ) = (-1)^\spr \sum_{i>k} (i-k) a_k = (-1)^{\spr }u_k
\end{equation*}
\end{lem}

\begin{proof}
Write 
\begin{equation*}
 C_k = \bigcup_{j<0}  D_{k,j}
\end{equation*}
 where \( D_{k,j} \) contains those 
elements of \( C_k \) of the form \( [x_i, j] \). From 
Lemma~\ref{lem:ck}, we see 
 that \( D_{k,j} \) is isomorphic to the complex \( \cfs (K ) \)
 truncated at level \( k-j  \) and translated down by \(2j\).
 By Proposition~\ref{prop:alex}, it contributes 
\begin{equation*}
 (-1)^\spr \sum _{i \geq k-j} a_j
\end{equation*}
to the Euler characteristic. Summing up these contributions for all 
\( j < 0 \) gives the desired result. 
\end{proof}

\noindent Putting these facts together, we obtain the following:

\begin{prop} For \( k \geq 0 \), there is an isomorphism of \( \Z [u] \) 
modules
\begin{align*}
H_{*} (C_k) \cong Q_k \oplus V_k 
\end{align*}
where \( Q_k \) and \(V_k \) are as in the statement of Theorem 1, except
that if \( \sigma \geq 0 \)  the grading
of \( u^{-h_k+1} \in V_k \) is \( \spr - 2. \)
\end{prop}

\begin{proof}
Suppose \( \sigma \geq 0 \). Then  for \( i < k-1 \) we have
 \( H_i (C_k) \cong \hfps _{i-2k} (K) = 0 \). Indeed, the latter group 
vanishes whenever \(i-2k < \spr \), and the assumption \( k\geq 0 \) 
implies \(i-2k < 0 \). 
Similarly, Corollary~\ref{cor:top} implies that for
 \( i > k-1 \), \( H_i(C_k) = \Z \) for 
\( i = \spr - 2, \spr - 4 \), and so on down to \(i = k \) or
\( i = k+1 \) (depending on the parity of \( \spr \)) and 
that the other homology groups are trivial.

 To get \(H_{k-1}(C_k) \),
we recall that it is free (by Lemma~\ref{lem:free}) and use
 Lemma~\ref{lem:chi} to compute its rank. We find that 
\(H_{k-1}(C_k) = \Z ^{|b_k|} \) if \( k \equiv \spr \mod 2 \)
and \( H_{k-1}(C_k) = \Z ^{|b_k|+1} \) otherwise. This proves the statement
at the level of groups. 
 
To check it at the level of modules, we note that the \(u \) 
action in the top half of \(C_k\) clearly agrees with the \(u \) 
action on \( \hfms
(K)\). In other words, it sends the generator in dimension \(\spr - 2i \)
to the generator in dimension \( \spr -2(i+1) \). In view of the structure
of \( H_*(C_k) \), it now suffices to show that if \( k < \spr \) and
 \( k-1 \equiv \spr
\mod 2 \), then \(u \) takes the generator of \( H_{k+1}(C_k ) \)
to a nonzero element of \(H_{k-1} (C_k) \). But this follows from the
fact that the \(u \) action commutes with the inclusion
\(  \hfms _i(K) \hookrightarrow H_i(C_k) \) for \(i \geq k-1 \). 

This proves the proposition in the case \( \sigma \geq 0 \). When
 \( \sigma \leq 0 \) the roles of the bottom and top halves are reversed,
but the proof is otherwise the same. 
\end{proof}

To prove Theorem 1, we use the long exact sequence:
$$\begin{CD}
 \cdots\qua H_i(C_k) @>>> \hfp _i (K, \spi _k) @>{\pi_*}>> \hfps _i (K)
@>>> H_{i-1}(C_k)\qua \cdots 
\end{CD}$$
We claim that the map \( \pi_* \) is always surjective.
This is clearly true when \(i \) is very large. But this implies the claim
for all \(i \), since \( \hfps _* (K) \cong \Zu \) and the maps in the 
long exact sequence respect the  \(u \) action. 

It is now a simple matter to check that the theorem holds at the level of 
groups. To verify the module structure, we consider the cases  
\( \sigma (K ) \geq 0 \) and \( \sigma (K ) \leq 0 \) separately.

\medskip

{\bf Case 1\qua \( \sigma (K ) \geq 0 \)}:\qua
From the short exact sequence of \( \Z [u] \) modules
\begin{equation*}
\label{eqss}
\begin{CD}
0 @>>>   H_*(C_k) @>>> HF_*^+ (K, \spi _k) @>>> \hfps _* (K) @>>> 0 
\end{CD}
\end{equation*}
we see that \( HF_*^+ (K, \spi _k) \) must have a  summand \( \Zu \)
which maps onto \( \hfps _* (K) \). We claim that the kernel of this map 
is isomorphic to \(V_k \). Indeed, the kernel can hardly be larger than
\(V_k\), since there is no such module available in \( H_*(C_k) \). 
To get other direction, we observe that \( CF_i^+ (K, \spi _k ) \) is
 identical to \( CF_i^{\infty } (K ) \) for \( i \geq k-1 \). Indeed,
the largest \(j\) for which \([x_i,j] \not \in \cfp (K, \spi _k) \) is
\(-1 \) if \(k \geq \epsilon _x (x_i) \) and
\(k- \epsilon _x (x_i) -1 \) if  \(k \leq \epsilon _x (x_i) \). In both
cases, \( \gr [x_i, j ] < k-1 \). It follows 
 \( HF_*^+ (K, \spi _k)\) has a \( \Zu \) summand which extends at least
 as far down in the grading as  \(V_k \). Together with the known
structure of \(\cfps (K) \), this implies that the kernel is at least as 
big as \(V_k\), thus proving the claim. It is now easy to see that 
to see that the \( \Z [u] \) module structure must be  as described. 

\medskip
{\bf Case 2\qua \( \sigma(K) \leq 0 \)}:\qua
As in the previous case,
\( HF_*^+ (K, \spi _k) \) must have a  summand \( \Zu \)
which maps onto \( \hfps _* (K) \). This time, however, it is easy
to show that this map has no kernel. Indeed, since \( \sigma \leq 0 \)
and \( k \geq 0 \), there are no elements in \( H_*(C_k) \) with grading
less than that of \( 1 \in \hfps _* (K) \cong \Zu \). Thus the short
exact sequence  splits as a sequence of \( \Z [u ] \) modules.
\qed

\

This finishes the proof of Theorem 1. Corollary~\ref{cor:hfhat}
 is an immediate  consequence  of the theorem and Lemma~\ref{lem:gysin},
while
  Corollary~\ref{cor:hfm} follows from the discussion of 
\( HF^{red} \) in section 4 of \cite{OS1}. 

\section{Using the exact triangle}
 
In this section, we describe how to calculate \( \hfp (K^n ) \) for
any \(n \). The results may be summarized as follows:

\begin{prop}
\label{prop:0surgery}
For \( k \neq 0 \), there is an isomorphism of 
 \( \Z [u ] \) modules
\begin{equation*}
 \hfp (K^0, \spi _k )  \cong Q_k \oplus V_k .
\end{equation*}
 Similarly, if we use the twisted coefficients of \cite{OS2} we have 
\begin{align*}
 \underline \hfp (K^0, \spi _0 ) &\cong (Q_0 \oplus V_0 ) \otimes \Z
[T,T^{-1}] \oplus \Zu. &
\end{align*}
\end{prop}

\begin{prop}
\label{prop:nsurgery}
For \( n  \in \Z^+ \) and \( \sigma (K ) \leq 0 \),  we have
\begin{equation*}
\hfp (K^n, \spi _k) \cong \Zu \oplus \bigoplus_{i \equiv k  (n)}
 Q_i \bigoplus _{i \equiv k (n)} V_i
\end{equation*}
while if  \( \sigma (K) \geq 0 \)
\begin{equation*}
\hfp (K^n, \spi _k) \cong \Zu \oplus \bigoplus_{i \equiv k  (n)} Q_i 
\bigoplus _{i_0 \neq i \equiv k  (n)} V_i
\end{equation*}
where \( i_0 \) is the representative of \( k \mod n \) with the smallest
absolute value. 
\end{prop}

{\bf Remarks}

\begin{enumerate}
\item The groups \(Q_k\) and \(V_k\) are described in the statement of 
Theorem 1. By convention,  \(Q_{-k} = Q_k \) and \(V_{-k} = V_k \).

\item For the moment, we consider all these groups to be \( \Z / 2 \) graded.
The absolute grading of \cite{OS3} can easily be worked out from the
exact triangle as well. 

\item We can use the same methods as in the large \(n \) case to compute
the groups \( \hfhat (K^n, \spi_k ) \) and \( \hfm (K^n, \spi _k ) \). To get
\( \hfp (K^{-n}, \spi _k ) \) we use the isomorphism 
\begin{equation*}
 \hfp (K^{-n}, \spi _k ) \cong (\hfm (\overline{K} ^n, \spi _k ))^* .
\end{equation*} 

\item These results are summarized in the following rule of thumb:
the Floer homology of the knot \(K \) behaves like the sum of a fixed
group \(Q \) and the Floer homology of the \((2,n) \) torus knot of the 
same signature. \(Q \) does not interact with \( \hfp (S^3) \)
in the exact triangle.
\end{enumerate}

To prove these two propositions, we use 
 the exact triangle of \cite{OS2}. Starting
with our knowledge of \( \hfp (K^n, \spi _k ) \) for \( n \gg 0 \), we 
use the exact triangle to calculate \( \hfp (K^0, \spi _k ) \). Then
we  apply the triangle again to get \( \hfp (K^n) \) for any \( n\). 

\subsection{Review of the exact triangle}

The exact triangle --- first conceived by Floer and now widened in scope
by Ozsv\'ath and Szab\'o --- has numerous incarnations. The one we want is 
described in Theorem 10.19 of \cite{OS2}.

\begin{thrm} For \( n \in \Z ^+ \), there are  exact triangles
\begin{equation*}
\xymatrix{ & \hfp (S^3) \ar[dl]& \\
\bigoplus_{i \equiv k (n) } \hfp (K^0, \spi _i ) \ar[rr] & & 
\hfp (K^n, \spi _k)  \ar[ul]}
\end{equation*}
and 
\begin{equation*}
\xymatrix{ & \hfp (S^3) \ar[dr]& \\
\bigoplus_{i \equiv k (n) } \hfp (K^0, \spi _i ) \ar[ur] & & 
\hfp (K^{-n}, \spi  _k ).  \ar[ll] }
\end{equation*}
All the maps in these triangles respect the \(\Z [u]\) module structure. 
\end{thrm}

{\bf Remarks}

\begin{enumerate}

\item Note that for \( n \geq 2g-1 \), the adjunction inequality
 implies that at most one term in the direct sum is  nontrivial.

\item Our notation for \( \spinc \) structures is as follows. In the language
of section 3.1, \( H^2(K^0) \) is generated by \( A_1 ^* \) (the Poincare 
dual of \(A_1 \)). On \(K^0 \), we let \( \spi _k \) be the \( \spinc \)
structure whose determinant line bundle has first Chern class \( 2kA_1^* \).
We have to check that for \( n \gg 0 \), the \( \spinc \) structure
\( \spi _k \) on \( K^0 \) really does induce \( \spi _k \) on \(K^n \). 
To see this, note that on \(K^n \), the difference \( \spi _{k+1} - \spi _k
= 2A_1^* \). It follows that there is a constant \(m \) such that
the \( \spinc \) structure \( \spi _k \) on \( K^0 \)  induces
\( \spi _{k+m} \) on \(K^n \). By adjunction,
\( \hfp (K^0, \spi _k ) \) is nontrivial only when \( -g < k < g \). 
On the other hand, we know that \( \hfp (K, \spi _k) \neq \hfp (S^3) \)
precisely for \( -g < k < g \), so we must have \( m = 0 \). 
(This proves that \( \spi _k \) and \( \spi _{-k} \) are conjugate
\( \spinc \) structures on \( K^n \), as mentioned in section 4.2.)

\item To use the second exact triangle, we need the following lemma.
\end{enumerate}
\begin{lem}
\label{lem:-surgery}
For \( n \gg 0 \), 
\begin{align*}
\hfp (K^{-n}, \spi _k ) \cong \begin{cases}
 \Zu \oplus Q_k \oplus V_k & \text{if} \  \sigma \geq 0 \\
  \Zu \oplus Q_k & \text{if} \  \sigma \leq 0 
\end{cases}
\end{align*}
\end{lem}

\begin{proof}
Observe that 
\begin{equation*}
\hfp (K^{-n}, \spi _k ) \cong HF _- ( {- K^{-n}}, \spi _k )
\cong (\hfm ( {\overline K}^{n} , \spi _k ))^*
\end{equation*}
where \( {\overline K } \) is the mirror image of \( K \). (In the last
equality, we have used the fact that 
\( \hfm ( {\overline K}^{n} , \spi _k) \)
is free over \( \Z \).) Since 
\( K \) and \( {\overline K } \) have the same Alexander polynomials
and \( \sigma ( {\overline K } ) = - \sigma (K) \), the lemma follows from
Corollary~\ref{cor:hfm}. 
\end{proof}

\subsection{Proof of Proposition~\ref{prop:0surgery}}
We use the exact triangle for \( n \gg 0 \).  
Consider the map \( h: \ \hfp (K^n, \spi _k)  \to \hfp (S^3 ) \cong \Zu \). 
For \( k \neq 0 \), \( \hfp (K^0, \spi _k ) \) is finitely generated, so
\( h \) must be onto when the grading is large. The only way this can
happen is if the \( \Zu \) summand in \( \hfp (K^n, \spi _k) \) maps
onto \( \hfp (S^3 ) \). In this case, \( \ker h \) is a torsion module
\( R = \Zu / u^{-m} \) plus the other summands in \( \hfp (K^n, \spi _k) \).

Suppose that \( \sigma (K) \geq 0 \), so that
 \( \hfp (K^n, \spi _k ) \cong \Zu \oplus Q_k  \). Then
\begin{equation*}
 \hfp (K^0, \spi _k ) \cong Q_k \oplus R .
\end{equation*} 
But by Theorem 9.1 of \cite{OS2} we have 
\( \chi (\hfp (K^0, \spi _k )) = u_k = \chi (H_*(C_k)) \). It follows that 
we must have \( R \cong V_k \). 
Similarly, if \( \sigma (K) \leq 0 \), we see that 
\( \hfp (K^0, \spi _k ) \cong Q_k \oplus V_k \oplus R \). In this case,
the Euler characteristic shows that the group \(R \) must be trivial. 

The case \( k = 0 \) requires some special attention, since 
\( \hfp (K^0, \spi _0 ) \) is not finitely generated. To get around this 
problem, we use the exact triangle for twisted coefficients (Theorem
10.21 of \cite{OS2}). In this case, we have a map 
\begin{equation*}
 h: \ \hfp (K^n, \spi _0)[T,T^{-1}] \to \hfp (S^3 )[T,T^{-1}]
 \cong \Zu [T,T^{-1}].
\end{equation*}
We consider the action of \(h\) on the \(\Zu \) summand of
 \(\hfp (K, \spi_0 )
\). Clearly \(h(u^{-k}) \) is nonzero --- otherwise \( \underline{\hfp }
(K^0, \spi _0 ) \) would contain a copy of\break \( \Zu [T,T^{-1}] \), which 
it does not. Thus we can write
\begin{equation*}
h(u^{-k}) = F(T) u^{-k+m} + (\text{lower order terms})
\end{equation*}
for some \(m \geq 0 \). We claim that after a change of basis in 
\(\hfp (S^3 )[T,T^{-1}] \), we have \(F(T) = T-1 \). This follows from the 
calculation of \(\hfi \) in \cite{OS2}. Indeed, the fact that 
\(\underline{\hfi} (K^0) \cong \Z [u,u^{-1}] \) implies that  (after a change
of basis) \(F \) is a polynomial of degree 
 1. On the other hand, to get \(h \) in the untwisted
triangle, we substitute \(T = 1 \) in the twisted map. Since 
\(\hfi (K^0) \cong H^*(T^2) \otimes \Z [u,u^{-1}] \), we must have 
\(F(1) = 0 \). 

Now that we  understand the action of \(h \), we see that
\begin{equation*}
 \underline{\hfp } (K^0, \spi _0) \cong (Q_0 \oplus R) \otimes \Z[T, T^{-1}] 
\oplus \Zu
\end{equation*}
where \(R \) is the  torsion module \( \Zu / u^{-m}\). 
Although it is not explicitly stated in \cite{OS2}, it is not difficult to 
follow through the proof of Theorem 9.1 there and see that it applies
to \( \underline{\hfp } (K^0, \spi _0) \) as well, so long as we take
the Euler characteristic of \( \underline{\hfp } (K^0, \spi _0) \)
as a \( \Z [T, T^{-1}] \) module. We
can now argue as in the case \( k \neq 0 \).
\qed 

\subsection{Proof of Proposition~\ref{prop:nsurgery}}
To compute \( \hfp (K^n, \spi _k ) \), we use the exact triangle
\begin{equation*}
\xymatrix{ & \hfp (S^3) \ar[dr]& \\
\bigoplus _{i \equiv k (n) } \hfp (K^0, \spi _i ) \ar[ur]^{\oplus g_i} & & 
\hfp (K^{-n}, \spi  _k ).  \ar[ll] }
\end{equation*}
The key point is that 
the maps \( g_{i} : \hfp (K^0, \spi _i ) \to \hfp (S^3) \) are independent
of \(n \). Indeed \( g_i \) is simply the map induced by the 
standard cobordism between \( S^3 \) and \( K^0 \).

\begin{lem}
\label{lem:g}
 For \( i \neq 0 \), the map \(g_i \) is trivial if \( \sigma \geq 0 
\). If \( \sigma \leq 0 \), \( g_i \) maps the summand
\(V_i \) onto the span of \( \{ 1, u^{-1}, u^{-2}, \ldots, u^{-h_i+1} \} \)
 in \( \hfp (S^3) \). If we use twisted coefficients, these statements
 are true for \( i=0 \) as well. With untwisted coefficients, 
 \( g_0 \) maps \( \hfp(K^0, \spi_0) \)  onto \( \hfp (S^3) \)
with kernel \(Q_0 \oplus V_0 \oplus \Zu \). 
\end{lem}

\begin{proof} Suppose \( i \neq 0 \) and take \( n \gg 0 \). 
 The usual argument shows that \( \hfp (S^3) \) must map 
onto the \( \Zu \) summand in \(\hfp (K^{-n}, \spi_k )\),
 with  a \( \Z [u] \) torsion module \(R \) as kernel.
For \( \sigma \geq 0 \), we have 
\( \hfp (K^{-n}, \spi_k ) \cong Q_k \oplus V_k \oplus \Zu \) by
 Lemma~\ref{lem:-surgery}. On the other hand, we
  know \( \hfp (K^0, \spi _i ) \cong Q_k \oplus V_k \),
so \(R \) must be trivial. 
Likewise, if \( \sigma \leq 0 \), we know
\( \hfp (K^{-n}, \spi_k ) \cong Q_k \oplus \Zu \), so we must have 
\( R \cong V_k \). 
For \( i = 0 \), the same arguments apply if we use twisted coefficients.

The statement for untwisted coefficients is just a rephrasing of the 
fact that the 
map \(\hfp (S^3) \to \hfp (K^{-n},\spi _0 ) \) is trivial.
\end{proof}

It is now elementary to compute \( \hfp (K^{-n}, \spi _k) \). For example,
if \( \sigma \geq 0 \) and \( k \not \equiv 0 \mod n \),
 all the \( g_i \) vanish, and we have 
\begin{equation*}
\hfp (K^{-n}, \spi _k) \cong \Zu \oplus \bigoplus_{i \equiv k  (n)} Q_i 
\bigoplus_{i \equiv k  (n) } V_i .
\end{equation*}
On the other hand, when \( \sigma \leq 0 \), the image of 
\(\oplus g_i \) is isomorphic to \( V_{k_0} \), where \( k_0 \) is the 
representative of \( k \mod n \) with the smallest absolute value. Thus
\begin{equation*}
\hfp (K^{-n}, \spi _k) \cong \Zu \oplus \bigoplus_{i\equiv k  (n) } Q_i 
\bigoplus_{k_0 \neq i \equiv k  (n)} V_i .
\end{equation*}
We leave it to the reader to check that these results also hold in the case
\( k \equiv 0 \mod n \). 

Finally, to get Proposition~\ref{prop:nsurgery}, we use the isomorphism
\begin{equation*}
 \hfp (K^{-n}, \spi _k) \cong (\hfm ({\overline K } ^n, \spi _k )) ^*
\end{equation*}
together with the usual exact sequence relating \( \hfp \) and \( \hfm \). 
\qed

\subsection{ Fr{\o}yshov's Invariant}
We conclude by justifying the statement about  Fr{\o}yshov's 
invariant given in the introduction. In \cite{OS3}, 
Ozsv\'ath and Szab\'o define an invariant \(d \) which is the analogue
(up to a factor of 2) of the {\it h}-invariant in their theory. 

To compute \( d (K^n) \), we use the fact that  the degree of the map
\( h\co \hfp (S^3 ) \to \hfp (K^{-n}, \spi _k )  \) 
 is determined by purely homological information,
namely the values of \(n\) and \(k\).  In particular, for \(n = 1\), the
degree of \(h \) is 0.

Consider the exact triangle for \( K^{-1} \). If we use twisted coefficients,
the  facts stated above imply that \( d (K^{-1} ) =  2 \dim \ker h \).
From Lemma~\ref{lem:g}, we know that
\begin{align*}
\dim \ker h =
\begin{cases} 
 0 & \text{if} \ \sigma \geq 0 \\
 |V_0| = \lceil |\spr | /2 \rceil & \text{if} \ \sigma  \leq 0 .
\end{cases}
\end{align*}
This gives us  \( d(K^{-1})\). To get
\( d(K^1) \), we use the easily proved fact that \(d({-Y} ) = 
- d(Y) \) and the identification \( { K^1 } = - {\overline K} ^{-1} \).

\section{Other Knots}
We conclude by discussing the extent to which the methods used in 
this paper extend to knots with more than two bridges. Of course there is 
no analogue of Schubert's theorem for these knots, and thus no explicit 
general form for their Heegard splitting. On the other hand, given a 
particular knot, it is straightforward to find a Heegard splitting of its 
complement. As in the case of two-bridge knots, it is convenient to
consider \(n\)-surgery on \(K\) for \(n \gg 0 \). In this case, we still have
a stable complex \( \cfs (K) \), with an Alexander grading which can be
computed by Fox calculus. As in section 3.3, it is straightforward
to determine \(\cfhat (K, \spi _k) \) from \( \cfs (K) \). In addition, we
should still be able to run the argument of section 4.3 to compute the 
grading of the generator of \( \cfs \). 

The calculation of section 4.4 fails, however, because for a general knot 
 the Alexander grading on \( \cfs (K) \) does not agree
with the Maslov grading. (For example, it is easy to see that the
two differ for non two-bridge torus knots.) We refer to those knots for 
which the two gradings coincide as {\it perfect} knots. We expect that the
methods  of this paper will extend to perfect knots in general,
although it is more difficult to  compute \( HF _a ^+ (K) \). 
We hope to discuss the issues described above, as well as some other
classes of perfect knots, in a future paper.

\Addresses\recd

\begin{thebibliography}

\bibitem{BZ} G. Burde and H. Zieschang, {\it Knots},
 de Gr{\"u}yter Studies in 
Mathematics, No. 5, Walter de Gr{\"u}yter and Co., Berlin, 1985. 

\bibitem{Fl} A. Floer, {\it Witten's complex and infinite dimensional Morse
theory}, J. Differential Geom. {\bf 30} (1989), 207--21. 

\bibitem{Fr} K. Fr{\o}yshov, Lectures on Floer homology at
 Harvard, Spring 2001. 

\bibitem{Mc} C. McMullen, {\it The Alexander polynomial of a 3-manifold 
and the Thurston norm on cohomology}, Ann. Scient. Ec. Norm. Sup., to 
appear.  

\bibitem{Mur} K. Murasugi, {\it Knot theory and its applications}, 
Birkhauser, Boston, 1996. 

\bibitem{OS1} P.S. Ozsv\'ath and Z. Szab\'o, {\it Holomorphic disks and 
topological invariants for homology three-spheres}, 
{\tt arXiv:math.SG/0101206}  

\bibitem{OS2} P.S. Ozsv\'ath and Z. Szab\'o, {\it Holomorphic disks and 
three-manifold invariants: properties and applications},
{\tt arXiv:math.SG/0105202}

\bibitem{OS3} P.S. Ozsv\'ath and Z. Szab\'o, {\it Absolutely graded Floer homologies 
and intersection forms for four-manifolds with boundary}, {\tt arXiv:math.SG/0110170}

\bibitem{Shu} H. Schubert, {\it Knoten mit zwei Br\"{u}cken}, Math. Z. {\bf 65}
(1956), 133-70.

\end{thebibliography}
\end{document}